\documentclass[11pt]{article}

\usepackage[english]{babel}

\usepackage{mathtools}    
\usepackage{amssymb}
\usepackage{amsthm}
\usepackage{bm}

\usepackage{fullpage}     

\usepackage{tikz}
\usepackage{tikz-cd}
\usetikzlibrary{positioning}
\usepackage{graphicx}

\usepackage{array}
\usepackage{multirow}

\usepackage{enumitem}

\usepackage{verbatim}
\usepackage[noadjust]{cite}

\usepackage{hyperref}
\hypersetup{hidelinks}

\renewcommand\arraystretch{1}

\setlist[itemize]{itemsep=0em}

\makeatletter
\renewcommand*\env@matrix[1][\arraystretch]{%
  \edef\arraystretch{#1}%
  \hskip -\arraycolsep
  \let\@ifnextchar\new@ifnextchar
  \array{*\c@MaxMatrixCols c}}
\makeatother

\theoremstyle{definition}
\newtheorem{definition}{Definition}[section]

\newtheorem{example}[definition]{Example}

\newtheorem{problem}[definition]{Problem}

\newtheorem{note}[definition]{Note}

\theoremstyle{plain}
\newtheorem{theorem}[definition]{Theorem}
\newtheorem{lemma}[definition]{Lemma}
\newtheorem{corollary}[definition]{Corollary}
\newtheorem{proposition}[definition]{Proposition}

\newcommand\scalemath[2]{\scalebox{#1}{\mbox{\ensuremath{\displaystyle #2}}}}

\DeclareMathOperator{\End}{End}		

\begin{document}

\title{\bf Some attractive bases for the Onsager Lie algebra}

\author {
Jae-Ho Lee\thanks{Department of Mathematics and Statistics, University of North Florida, Jacksonville, FL 32224, USA.}
}
\date{}

\maketitle
\newcommand\blfootnote[1]{%
\begingroup
\renewcommand\thefootnote{}\footnote{#1}%
\addtocounter{footnote}{-1}%
\endgroup}
\blfootnote{E-mail address:  \texttt{jaeho.lee@unf.edu}}

\begin{abstract}
\noindent
The Onsager Lie algebra $\mathcal{O}$ is defined by generators $A$, $B$ subject to the Dolan-Grady relations $[A,[A,[A,B]]] = 4[A,B]$ and $[B,[B,[B,A]]] = 4[B,A]$. It is known that $\mathcal{O}$ is infinite-dimensional.
In this paper, we consider three bases for $\mathcal{O}$.
The first basis was recently introduced by the author, and the second basis is a variation of the first.
The third basis was introduced by Baseilhac and Cramp\'e in 2019.
By construction, the elements in each basis are recursively obtained from $A$, $B$. 
Using the $\mathfrak{sl}_2$ loop algebra, we express every basis element in closed form.
We then describe how these three bases are related to each other.
Next, we consider a finite-dimensional irreducible $\mathcal{O}$-module $V$.
Associated with $V$ are six well-known decompositions: the eigenspace decomposition of $A$ and the eigenspace decomposition of $B$, together with four other decompositions of $V$ said to be split.
We show that on each decomposition, every element of the first two bases acts in a block diagonal, block upper bidiagonal, block lower bidiagonal, or block tridiagonal fashion.
We illustrate these actions using a tetrahedron diagram.
We also determine how the elements of the third basis act on the six decompositions.
To obtain our results about finite-dimensional irreducible $\mathcal{O}$-modules, we use the theory of tridiagonal pairs.


\end{abstract}
\section{Introduction}
This paper concerns the Onsager Lie algebra and its finite-dimensional irreducible modules.
In his seminal paper on the two-dimensional Ising model in zero magnetic field, Onsager introduced an infinite-dimensional Lie algebra, now called the Onsager Lie algebra \cite{1944Ons}.
In their study of conserved quantities in integrable systems, Dolan and Grady introduced a pair of relations involving two operators \cite{1982DG}.
Perk later used these relations to obtain a presentation of the Onsager Lie algebra \cite{1989Perk}.
In Perk's presentation, the Onsager Lie algebra $\mathcal{O}$ is defined by generators $A$, $B$ subject to the relations:
\begin{align*}
	[A,[A,[A,B]]] & = 4[A,B], \\
	[B,[B,[B,A]]] & = 4[B,A].
\end{align*}
We call $A$, $B$ the standard generators of $\mathcal{O}$.
We now mention some results about $\mathcal{O}$.
Roan realized $\mathcal{O}$ as a Lie subalgebra of the $\mathfrak{sl}_2$ loop algebra \cite{1991Roan}.
Baseilhac, Belliard, and Cramp\'e obtained a presentation of $\mathcal{O}$ in the Faddeev-Reshetikhin-Takhtajan (FRT) framework \cite{2018BBC}.
More recently, we obtained four bases for $\mathcal{O}$ that are related by a $\mathbb{Z}_2\times\mathbb{Z}_2$ action \cite{2026LeeJOAC}.

\smallskip
The Onsager Lie algebra $\mathcal{O}$ is closely related to tridiagonal pairs.
Roughly speaking, a tridiagonal pair (or TD pair) consists of two diagonalizable linear maps on a nonzero finite-dimensional vector space, each of which acts in a block tridiagonal fashion on the eigenspaces of the other \cite{ITT2001}.
Hartwig showed that the standard generators of $\mathcal{O}$ act as a TD pair on every finite-dimensional irreducible $\mathcal{O}$-module \cite{2007Har}.

\smallskip
In the present paper, we consider three bases for $\mathcal{O}$:
\begin{alignat}{4}
	& A_k,     \qquad && B_k,     \qquad && \psi_{k+1} \qquad && (k \in \mathbb{N}); \label{eq:intro first basis}\\
	& A_k,     \qquad && B_k,     \qquad && \tilde\psi_{k+1} \qquad && (k \in \mathbb{N}); \label{eq:intro second basis}\\
	& W_{-k},  \qquad && W_{k+1}, \qquad && \tilde G_{k+1} \qquad && (k \in \mathbb{N}). \label{eq:intro third basis}
\end{alignat}
The basis \eqref{eq:intro first basis} was introduced by the author in  \cite{2026LeeJOAC}.
The basis \eqref{eq:intro second basis} is introduced in this paper as a variation of the basis \eqref{eq:intro first basis}.
The bases \eqref{eq:intro first basis} and \eqref{eq:intro second basis} are related by the equation
\begin{equation*}
	\tilde\psi_{k+1}=A_k+B_k-\psi_{k+1} \qquad (k\in\mathbb{N}).
\end{equation*}
The basis \eqref{eq:intro third basis} was introduced by Baseilhac and Cramp\'e in their work on the FRT presentation of $\mathcal{O}$ \cite{2019BC}.
By construction, the elements of each basis in \eqref{eq:intro first basis}--\eqref{eq:intro third basis} are recursively obtained from the standard generators $A$, $B$.

\smallskip
We now summarize our main results.
We have two types of results.
The first type is about the Lie algebra $\mathcal{O}$, and the second type is about finite-dimensional irreducible $\mathcal{O}$-modules.
The main results of the first type are as follows.
We use Roan's embedding of $\mathcal{O}$ into the $\mathfrak{sl}_2$ loop algebra to express the elements in \eqref{eq:intro first basis}--\eqref{eq:intro third basis} in closed form as $2\times 2$ matrices.
We then use these matrix expressions to determine the transition matrices between the bases \eqref{eq:intro first basis}--\eqref{eq:intro third basis}.
The main results of the second type are as follows.
Let $V$ denote a finite-dimensional irreducible $\mathcal{O}$-module.
There are six well-known decompositions of $V$: the eigenspace decomposition of the $\mathcal{O}$-generator $A$ and the eigenspace decomposition of the $\mathcal{O}$-generator $B$, together with four other decompositions of $V$ that are called split.
We determine how each element of the bases \eqref{eq:intro first basis}--\eqref{eq:intro third basis} acts on these six decompositions of $V$.
We show that on each of the six decompositions, every element of the bases \eqref{eq:intro first basis} and \eqref{eq:intro second basis} acts in a block diagonal, block upper bidiagonal, block lower bidiagonal, or block tridiagonal fashion.
We illustrate these actions using a tetrahedron diagram.
The action of the basis elements in \eqref{eq:intro third basis} is described as follows.
We show that for $k \in \mathbb{N}$ the action of $W_{-k}$ is similar to that of $A_k$, the action of $W_{k+1}$ is similar to that of $B_k$, and the action of $\tilde G_{k+1}$ is block tridiagonal on each of the six decompositions.
To obtain our results about finite-dimensional irreducible $\mathcal{O}$-modules, we use the theory of TD pairs.

\smallskip
Our study of finite-dimensional irreducible $\mathcal{O}$-modules is motivated by the work of Terwilliger \cite{2023TerJCTA}.
In \cite[Sections~8,9]{2023TerJCTA}, he considered the positive part $U_q^+$ of the quantum affine algebra $U_q(\widehat{\mathfrak{sl}}_2)$.
He described how the alternating elements of $U_q^+$ act on six decompositions of a finite-dimensional irreducible $U_q^+$-module of non-nil type \cite{IT2007JAA}.
The results of the second type in the present paper can be viewed as a $q=1$ analogue of the results about $U_q^+$ in \cite{2023TerJCTA}.

\smallskip
This paper is organized as follows.
In Section~\ref{Sec:prelim}, we give some preliminaries.
In Section~\ref{Sec:three-bases}, we present the three bases \eqref{eq:intro first basis}--\eqref{eq:intro third basis} and obtain some relations in $\mathcal{O}$.
In Section~\ref{Sec:matrix-rep}, we represent the basis elements in \eqref{eq:intro first basis}--\eqref{eq:intro third basis} as $2\times 2$ matrices in the $\mathfrak{sl}_2$ loop algebra.
In Section~\ref{Sec:TM}, we determine the transition matrices between the three bases \eqref{eq:intro first basis}--\eqref{eq:intro third basis}.
In Sections~\ref{Sec:TD_pair}--\ref{Sec:decomp-diag}, we review TD pairs and TD systems, along with decompositions and the tetrahedron diagram.
In Sections~\ref{sec:TDsystemO-mod} and~\ref{Sec:action-on-V}, we apply the theory of TD pairs to a finite-dimensional irreducible $\mathcal{O}$-module $V$. 
We determine how the elements in \eqref{eq:intro first basis} and \eqref{eq:intro second basis} act on the six decompositions of $V$.
We illustrate these actions using the tetrahedron diagram.
In Section~\ref{Sec:future-research}, we give three open problems.

\section{Preliminaries}\label{Sec:prelim}
In this section, we fix some notation and recall a few basic facts that will be used throughout the paper.
Let $\mathbb{F}$ denote a field of characteristic zero.
Recall the set of natural numbers $\mathbb{N} = \{0,1,2,\ldots \}$.
Every vector space and tensor product in this paper is over $\mathbb{F}$. 
Every algebra and Lie algebra in this paper is over $\mathbb{F}$. 
Every algebra mentioned without the ``Lie'' prefix is associative and has a multiplicative identity.
Let $V$ denote a nonzero finite-dimensional vector space.
Let $\operatorname{End}(V)$ denote the algebra of all $\mathbb{F}$-linear maps from $V$ to $V$. 
Let $I$ denote the identity element in $\operatorname{End}(V)$.
An element $A \in \operatorname{End}(V)$ is said to be \emph{diagonalizable} whenever $V$ is spanned by the eigenspaces of $A$. 
Assume that $A$ is diagonalizable. 
Let $\{V_i\}_{i=0}^d$ denote an ordering of the eigenspaces of $A$.
For $0 \leq i \leq d$, let $\theta_i$ denote the eigenvalue of $A$ associated with $V_i$. 
For $0 \leq i \leq d$, define $E_i \in \operatorname{End}(V)$ by
$$
	(E_i - I )V_i = 0, \qquad 
	E_iV_j  = 0 \qquad(0 \leq j \leq d, \ j \neq i).
$$
We call $E_i$ the \emph{primitive idempotent} of $A$ associated with $V_i$ (or $\theta_i$).
By linear algebra, the following hold:
\begin{align*}
	& E_iE_j  = \delta_{i,j}E_i \quad (0 \leq i, j \leq d), \qquad
	I  = \sum_{i=0}^d E_i, \qquad
	A  = \sum^d_{i=0}\theta_iE_i, \\ 
	& AE_i  = E_iA = \theta_iE_i \quad (0 \leq i \leq d), \qquad
	V_i  = E_iV \quad (0 \leq i \leq d).
\end{align*}
Moreover,
\begin{equation*}
	E_i = \prod_{\substack{0 \leq j \leq d \\ j \neq i}} \frac{A - \theta_j I}{\theta_i - \theta_j} \qquad (0 \leq i \leq d).
\end{equation*} 
We note that each of $\{E_i\}_{i=0}^d$ and $\{A^i\}_{i=0}^d$ is a basis for the subalgebra of $\operatorname{End}(V)$ generated by $A$.

\section{Three bases for the Onsager Lie algebra $\mathcal{O}$}\label{Sec:three-bases}

In this section, we discuss three bases for the Onsager Lie algebra $\mathcal{O}$.
We begin by recalling the definition of $\mathcal{O}$.

\begin{definition}[\cite{1944Ons,1989Perk}]\label{def:Onsager}
Let $\mathcal{O}$ denote the Lie algebra with generators $A,B$ and relations
\begin{align}
	[A,[A,[A,B]]] & =4[A,B], \label{OnsagerDG(1)}\\
	[B,[B,[B,A]]] & =4[B,A]. \label{OnsagerDG(2)}
\end{align}
We call $\mathcal{O}$ the \emph{Onsager Lie algebra}.
We call $A,B$ the \emph{standard generators} of $\mathcal{O}$.
We refer to $A$ (resp. $B$) as the \emph{first standard generator} (resp. \emph{second standard generator}).
We call the relations \eqref{OnsagerDG(1)} and \eqref{OnsagerDG(2)} the \emph{Dolan-Grady relations}.
\end{definition}

We first recall the following basis for $\mathcal{O}$.

\begin{lemma}[{\cite[Theorem~6.3]{2026LeeJOAC}}]\label{O basis}
There exists a basis
\begin{equation}\label{eq:1st basis O}
	A_k, \qquad B_k, \qquad \psi_{k+1}, \qquad \qquad k \in \mathbb{N}
\end{equation}
for $\mathcal{O}$ such that $A_0 = A$, $B_0 = B$, and for $k, \ell \in \mathbb{N}$,
\begin{align}
	[A_k, A_\ell] & = 0, \qquad
	[B_k, B_\ell] = 0, \qquad
	[\psi_{k+1}, \psi_{\ell+1}] = 0, \label{thm0312:eq[A,B,psi]}\\[0.2em]
	[\psi_{k+1}, A_\ell] & = 2\psi_{k+\ell+1} + 2A_{k+\ell+1}, \label{thm0312:eq[psi,A]}\\[0.2em]
	[B_k, \psi_{\ell+1}] & = 2B_{k+\ell+1} + 2\psi_{k+\ell+1}, \label{thm0312:eq[B,psi]}\\[0.2em]
	[A_k, B_\ell] & = 2A_{k+\ell} + 2B_{k+\ell} - 4\psi_{k+\ell+1}. \label{thm0312:eq[A,B]}
\end{align}
\end{lemma}

The next lemma gives a construction of the basis \eqref{eq:1st basis O} from $A$ and $B$.

\begin{lemma}[{\cite[Corollary~6.5]{2026LeeJOAC}}]\label{lem:1st-rec}
The basis \eqref{eq:1st basis O} is recursively obtained from $A$ and $B$ in the order 
$$
	A_0, \quad B_0, \quad \psi_1, \quad 
	A_1, \quad B_1, \quad \psi_2, \quad 
	A_2, \quad B_2, \quad \psi_3, \quad \ldots
$$
using the following equations: for $k \in \mathbb{N}$,
\begin{align}
	\psi_{k+1} & = \frac{A_{k}}{2} + \frac{B_{k}}{2} - \frac{[A_{k}, B]}{4} = \frac{A_{k}}{2} + \frac{B_{k}}{2} - \frac{[A, B_k]}{4}, \label{cor:0312 psi}\\
	A_{k+1} & = \frac{[\psi_{k+1}, A]}{2} - \psi_{k+1}, \label{cor:0312 A}\\
	B_{k+1} & = \frac{[B, \psi_{k+1}]}{2} - \psi_{k+1}. \label{cor:0312 B}
\end{align}
\end{lemma}

\begin{example}
Using Lemma~\ref{lem:1st-rec}, we express $A_k, B_k, \psi_{k+1}$ $(0 \leq k \leq 2)$ in terms of $A$ and $B$.
\begin{align*}
	A_0 & = A,\\
	B_0 & = B, \\
	\psi_{1} & = \frac{A}{2} + \frac{B}{2} - \frac{[A,B]}{4}, \\[2pt]
	A_{1} & = -\frac{A}{2} -\frac{B}{2} - \frac{[A,[B,A]]}{8}, \\[2pt]
	B_{1} & = -\frac{A}{2} -\frac{B}{2} - \frac{[B,[A,B]]}{8}, \\[2pt]
	\psi_{2} & = -\frac{A}{2} - \frac{B}{2} - \frac{[B,A]}{8} - \frac{[A, [B, A]]}{16} - \frac{[B,[A,B]]}{16} - \frac{[B, [A, [B, A]]]}{32}, \\[2pt]
	A_{2} & = \frac{A}{2} + \frac{B}{2} + \frac{[A, [B, A]]}{8} + \frac{[B, [A, B]]}{16} + \frac{[A, [B, [A, [B, A]]]]}{64}, \\[2pt]
	B_{2} & = \frac{A}{2} + \frac{B}{2} + \frac{[B, [A, B]]}{8} + \frac{[A, [B, A]]}{16} + \frac{[B, [A, [B, [A, B]]]]}{64},\\[2pt]
	\psi_{3} & = \frac{A}{2} + \frac{B}{2} + \frac{[B,A]}{16} + \frac{3[A,[B,A]]}{32} + \frac{3[B,[A,B]]}{32} + \frac{[B, [A,[B,A]]]}{32} \\[2pt]
	& \qquad + \frac{[A,[B,[A,[B,A]]]]}{128}  + \frac{[B,[A,[B,[A,B]]]]}{128} + \frac{[B,[A,[B,[A,[B,A]]]]]}{256}.
\end{align*}
\end{example}

We now introduce a second basis for $\mathcal{O}$.
As we will see, the second basis closely parallels the first.

\begin{definition}\label{def:tilde psi_i}
For $k \in \mathbb{N}$, define $\tilde \psi_{k+1} \in \mathcal{O}$ by
\begin{equation}\label{eq:tilde psi_i}
	\tilde\psi_{k+1} = A_k + B_k - \psi_{k+1}.
\end{equation}
\end{definition}

\begin{lemma}\label{O basis-2}
With reference to Definition \ref{def:tilde psi_i}, the elements
\begin{equation}\label{eq:2nd basis O}
	A_k, \qquad B_k, \qquad \tilde\psi_{k+1}, \qquad \qquad k \in \mathbb{N}
\end{equation}
form a basis for $\mathcal{O}$.
Moreover, for $k, \ell \in \mathbb{N}$,
\begin{align}
	[A_k, A_\ell] & = 0, \qquad
	[B_k, B_\ell] = 0, \qquad
	[\tilde\psi_{k+1}, \tilde\psi_{\ell+1}] = 0, \label{eq[A,B,t-psi]}\\[0.2em]
	[A_k, \tilde\psi_{\ell+1}] & = 2A_{k+\ell+1} + 2\tilde\psi_{k+\ell+1}, \label{eq[t-psi,A]}\\[0.2em]
	[\tilde\psi_{k+1}, B_\ell] & = 2\tilde\psi_{k+\ell+1} + 2B_{k+\ell+1} , \label{eq[B,t-psi]}\\[0.2em]
	[B_k, A_\ell] & = 2B_{k+\ell} + 2A_{k+\ell} - 4\tilde\psi_{k+\ell+1}. \label{eq[B,A], t-psi}
\end{align}
\end{lemma}
\begin{proof}
By Lemma \ref{O basis} and Definition \ref{def:tilde psi_i}.
\end{proof}

\begin{lemma}\label{lem:2nd-rec}
The basis \eqref{eq:2nd basis O} is recursively obtained from $A$ and $B$ in the order
$$
	A_0, \quad B_0, \quad \tilde\psi_1, \quad 
	A_1, \quad B_1, \quad \tilde\psi_2, \quad 
	A_2, \quad B_2, \quad \tilde\psi_3, \quad \ldots
$$
using the following equations: for $k \in \mathbb{N}$,
\begin{align}
	\tilde\psi_{k+1} & = \frac{A_k}{2} + \frac{B_k}{2} + \frac{[A_k, B]}{4} = \frac{A_k}{2} + \frac{B_k}{2} + \frac{[A, B_k]}{4}, \label{eq:rec tpsi}\\
	A_{k+1} & = \frac{[A, \tilde\psi_{k+1}]}{2} - \tilde\psi_{k+1}, \notag \\
	B_{k+1} & = \frac{[\tilde\psi_{k+1}, B]}{2} - \tilde\psi_{k+1}. \notag
\end{align}
\end{lemma}
\begin{proof}
By \eqref{eq[t-psi,A]}--\eqref{eq[B,A], t-psi}.
\end{proof}

\begin{example}
Using Lemma~\ref{lem:2nd-rec}, we express $A_k, B_k, \tilde\psi_{k+1}$ $(0 \leq k \leq 2)$ in terms of $A$ and $B$.
\begin{align*}
	A_0 & = A, \\[2pt]
	B_0 & = B, \\[2pt]
	\tilde\psi_1 & = \frac{A}{2} + \frac{B}{2} + \frac{[A,B]}{4}, \\[2pt]
	A_1 & = - \frac{A}{2} - \frac{B}{2} - \frac{[A, [B, A]]}{8}, \\[2pt]
	B_1 & = - \frac{A}{2} - \frac{B}{2} - \frac{[B, [A, B]]}{8}, \\[2pt]
	\tilde\psi_2 & = - \frac{A}{2} - \frac{B}{2} - \frac{[A,B]}{8} - \frac{[A, [B, A]]}{16}  - \frac{[B, [A, B]]}{16}  - \frac{[A, [B, [A, B]]]}{32} , \\[2pt]
	A_2 & = \frac{A}{2} + \frac{B}{2} + \frac{[A, [B, A]]}{8} + \frac{[B, [A,B]]}{16} + \frac{[A, [B, [A, [B, A]]]]}{64}, \\[2pt]
	B_2 & = \frac{A}{2} + \frac{B}{2} + \frac{[B, [A, B]]}{8} + \frac{[A, [B,A]]}{16} + \frac{[B, [A, [B, [A, B]]]]}{64}, \\[2pt]
	\tilde\psi_3 & = \frac{A}{2} + \frac{B}{2} + \frac{[A,B]}{16} + \frac{3[A, [B, A]]}{32}  + \frac{3[B, [A, B]]}{32} + \frac{[A, [B, [A, B]]]}{32}  \\[2pt]
	& \qquad + \frac{[A, [B, [A, [B, A]]]]}{128}  + \frac{[B, [A, [B, [A, B]]]]}{128}  + \frac{[A, [B, [A, [B, [A, B]]]]]}{256}.
\end{align*}
\end{example}

So far, we have given two bases for $\mathcal{O}$.
We have a comment on a relation between these two bases.

\begin{lemma}\label{lem:PT5.7}
For $k, \ell \in \mathbb{N}$ the following hold in $\mathcal{O}$:
\begin{align}
	[A_k, B_\ell] & = 2(\tilde\psi_{k+\ell+1} - \psi_{k+\ell+1}), \label{eq:[Ai,Bj]} \\
	[\psi_{k+1}, \tilde\psi_{\ell+1}] & = 2(A_{k+\ell+1} - B_{k+\ell+1}). \label{eq:[psii,tpsij]}
\end{align}
\end{lemma}
\begin{proof}
Use \eqref{thm0312:eq[A,B]} and \eqref{eq:tilde psi_i} to get \eqref{eq:[Ai,Bj]}.
Evaluate $[\psi_{k+1}, \tilde\psi_{\ell+1}]$ using \eqref{eq:tilde psi_i}, and simplify the result using \eqref{thm0312:eq[psi,A]} and \eqref{thm0312:eq[B,psi]} to get \eqref{eq:[psii,tpsij]}.
\end{proof}

The relations in the next lemma are motivated by analogous relations for the alternating central extension of $\mathcal{O}$; see \cite[Lemma~4.24]{2023TerCAlg}.

\begin{lemma}\label{lem:rels in O}
For $k, \ell \in \mathbb{N}$ the following relations hold in $\mathcal{O}$:
\begin{align*}
	& [A_k, A_\ell] = 0, \qquad [B_k, B_\ell] = 0, \\
	& [A_k, B_\ell] + [B_k, A_\ell] = 0, \\
	& [A_k, \psi_{\ell+1}] + [\psi_{k+1}, A_\ell] = 0, \\
	& [A_k, \tilde\psi_{\ell+1}] + [\tilde\psi_{k+1}, A_\ell] = 0, \\
	& [B_k, \psi_{\ell+1}] + [\psi_{k+1}, B_\ell] = 0, \\
	& [B_k, \tilde\psi_{\ell+1}] + [\tilde\psi_{k+1}, B_\ell] = 0, \\
	& [\psi_{k+1}, \psi_{\ell+1}] = 0, \qquad [\tilde\psi_{k+1}, \tilde\psi_{\ell+1}] = 0, \\
	& [\tilde\psi_{k+1}, \psi_{\ell+1}] + [\psi_{k+1}, \tilde\psi_{\ell+1}] = 0.
\end{align*}
\end{lemma}
\begin{proof}
Use Lemmas \ref{O basis}, \ref{O basis-2} along with Definition \ref{def:tilde psi_i}.
\end{proof}

The following relations will be used in Section~\ref{sec:TDsystemO-mod}.
\begin{lemma}\label{lem:some rels A,B}
For $k \in \mathbb{N}$, the following {\rm(i)}--{\rm(iv)} hold in $\mathcal{O}$.
\begin{itemize}
	\item[\rm(i)] $A_k$ satisfies
	$$
		[A, A_k] = 0, \qquad \qquad [B, [B, [B, A_k]]] - 4[B, A_k] = 0.
	$$
	\item[\rm(ii)] $B_k$ satisfies
	$$
		[A, [A, [A, B_k]]] - 4[A, B_k] = 0, \qquad \qquad [B, B_k] = 0.
	$$
	\item[\rm(iii)] $\psi_{k+1}$ satisfies
	$$
		[A, [A, \psi_{k+1}]] + 2[A, \psi_{k+1}] = 0, \qquad \qquad 
		[B, [B, \psi_{k+1}]] - 2[B, \psi_{k+1}] = 0.
	$$
	\item[\rm(iv)] $\tilde\psi_{k+1}$ satisfies
	$$
		[A, [A, \tilde\psi_{k+1}]] - 2[A, \tilde\psi_{k+1}] = 0, \qquad \qquad 
		[B, [B, \tilde\psi_{k+1}]] + 2[B, \tilde\psi_{k+1}] = 0.
	$$
\end{itemize}
\end{lemma}

\begin{proof}
(i) The first equation follows immediately from \eqref{thm0312:eq[A,B,psi]}.
For the second equation, evaluate $[B, [B, [B, A_k]]]$ using \eqref{eq[A,B,t-psi]}, \eqref{eq[B,t-psi]}, and \eqref{eq[B,A], t-psi}. 
Simplify the result to obtain $4[B, A_k]$. \\
(ii) This is similar to (i).\\
(iii) For the first equation, evaluate $[A, [A, \psi_{k+1}]]$ using \eqref{thm0312:eq[A,B,psi]} and \eqref{thm0312:eq[psi,A]}.
Simplify the result to obtain $-2[A, \psi_{k+1}]$.  
The second equation is obtained similarly.\\
(iv) For the first equation, evaluate $[A, [A, \tilde \psi_{k+1}]]$ using \eqref{eq[A,B,t-psi]} and \eqref{eq[t-psi,A]}.
Simplify the result to obtain $2[A, \tilde\psi_{k+1}]$.
The second equation is obtained similarly.
\end{proof}

We have discussed the two bases for $\mathcal{O}$ given in \eqref{eq:1st basis O} and \eqref{eq:2nd basis O}.
We now present a third basis.
This basis was introduced by Baseilhac and Cramp\'e \cite{2019BC} in order to study the structure of $\mathcal{O}$ and its connections with integrable systems.
We will use the normalization of this basis given in \cite{2023TerCAlg}.

\begin{lemma}[{\cite[Definition~4.1,~Theorem~2]{2019BC}, \cite[Lemma~3.8]{2023TerCAlg}}]
There exists a basis 
\begin{equation}\label{eq:3rd basis O}
	W_{-k}, \qquad W_{k+1}, \qquad \tilde{G}_{k+1}, \qquad \qquad k \in \mathbb{N}
\end{equation}
for $\mathcal{O}$ such that $W_0 = A$, $W_1=B$, and for $k, \ell \in \mathbb{N}$,
\begin{align}
	& [W_{-k}, W_{-\ell}] = 0 , \qquad
	[W_{k+1}, W_{\ell+1}] = 0, \qquad
	[\tilde{G}_{k+1}, \tilde{G}_{\ell+1}] = 0, \notag \\[0.2em]
	& [W_{-k}, W_{\ell+1}] = \tilde{G}_{k+\ell+1}, \label{eq:baO(1)}\\[0.2em]
	& [\tilde{G}_{k+1}, W_{-\ell}] = 4W_{-k-\ell-1} - 4W_{k+\ell+1}, \label{eq:baO(2)}\\[0.2em]
	& [W_{k+1}, \tilde{G}_{\ell+1}] = 4W_{k+\ell+2} - 4W_{-k-\ell}.\label{eq:baO(3)}
\end{align}
\end{lemma}

\begin{lemma}\label{lem:3rd-rec}
The basis \eqref{eq:3rd basis O} is recursively obtained from $A$ and $B$ in the order
$$
	W_0, \quad W_1, \quad 
	\tilde{G}_1, \quad W_{-1}, \quad W_2, 
	\quad \tilde{G}_2, \quad W_{-2}, \quad W_3, 
	\quad \tilde{G}_3, \quad W_{-3}, \quad W_4,  \quad \ldots
$$
using the following equations: for $k \in \mathbb{N}$,
\begin{align*}
	\tilde{G}_{k+1} & = [W_0, W_{k+1}] = [W_{-k}, W_1], \\[2pt]
	W_{-k-1} & = W_{k+1} + \frac{[\tilde{G}_{k+1}, W_0]}{4}, \\[2pt]
	W_{k+2} & = W_{-k} + \frac{[W_1, \tilde{G}_{k+1}]}{4}.
\end{align*}
\end{lemma}
\begin{proof}
By \eqref{eq:baO(1)}--\eqref{eq:baO(3)}.
\end{proof}

\begin{example}
Using Lemma~\ref{lem:3rd-rec}, we express $W_{-k}$, $W_{k+1}$, $\tilde{G}_{k+1}$ $(0 \leq k \leq 2)$ in terms of $A$, $B$.
\begin{align*}
	W_0 & = A, \\[2pt]
	W_{1} & = B, \\[2pt]
	\tilde{G}_1 & = [A, B], \\[2pt]
	W_{-1} & = B + \frac{[A, [B, A]]}{4}, \\[2pt]
	W_2 & = A + \frac{[B, [A, B]]}{4}, \\[2pt]
	\tilde{G}_2 & = \frac{[A, [B, [A, B]]]}{4}  = - \frac{[B, [A, [B, A]]]}{4}, \\[2pt]
	W_{-2} & = A + \frac{[B, [A, B]]}{4}  + \frac{[A, [B, [A, [B, A]]]]}{16} , \\[2pt]
	W_3 & = B + \frac{[A, [B, A]]}{4}  + \frac{[B, [A, [B, [A, B]]]]}{16} , \\[2pt]
	\tilde{G}_3 & = \frac{[A, [B, [A, [B, [A, B]]]]]}{16}  = - \frac{[B, [A, [B, [A, [B, A]]]]]}{16}.
\end{align*}
\end{example}

In this section, we gave three bases for $\mathcal{O}$.
These bases are displayed in \eqref{eq:1st basis O}, \eqref{eq:2nd basis O}, \eqref{eq:3rd basis O}.
In Section \ref{Sec:TM} we will describe how these three bases are related to each other.

\section{Matrix representations of $\mathcal{O}$}\label{Sec:matrix-rep}
We continue to discuss the three bases for $\mathcal{O}$ given in \eqref{eq:1st basis O}, \eqref{eq:2nd basis O}, \eqref{eq:3rd basis O}.
In this section, we use Roan's embedding of $\mathcal{O}$ into the $\mathfrak{sl}_2$ loop algebra to describe the images of the elements in the three bases.
To this end, we first describe the $\mathfrak{sl}_2$ loop algebra.

\smallskip
Recall that $\mathfrak{sl}_2$ is the vector space consisting of all $2\times 2$ matrices with trace zero, equipped with the Lie bracket $[u,v]=uv-vu$ for $u,v\in\mathfrak{sl}_2$.
Let $t$ denote an indeterminate. 
Let $L=L(\mathfrak{sl}_2)$ denote the Lie algebra consisting of the vector space $\mathfrak{sl}_2 \otimes \mathbb{F}[t,t^{-1}]$ and the Lie bracket
$$
	[u \otimes \phi, v \otimes \varphi] = [u,v] \otimes \phi \varphi, \qquad \qquad u,v \in \mathfrak{sl}_2, \qquad \phi, \varphi \in \mathbb{F}[t,t^{-1}].
$$
We call $L$ the \emph{$\mathfrak{sl}_2$ loop algebra}.
We regard each element $m\in L$ as a trace-zero $2\times2$ matrix with entries in $\mathbb{F}[t, t^{-1}]$:
$$
	m = \begin{pmatrix}
	a & b \\
	c & d
	\end{pmatrix},
	\qquad \quad
	a, b, c, d \in \mathbb{F}[t,t^{-1}],
	\quad a+d = 0.
$$
By \cite[Proposition~1]{1991Roan}, there exists an injective Lie algebra homomorphism $\xi: \mathcal{O} \to L$ that sends
\begin{equation}\label{eq:image of A, B}
	A \ \longmapsto \ \begin{pmatrix} 0 & 1 \\ 1 & 0 \end{pmatrix}, \qquad
	B \ \longmapsto \ \begin{pmatrix} 0 & t \\ t^{-1} & 0 \end{pmatrix}.
\end{equation}
Recall the elements of $\mathcal{O}$ from \eqref{eq:1st basis O} and \eqref{eq:2nd basis O}:
\begin{equation*}
	A_k, \qquad B_k, \qquad \psi_{k+1}, \qquad \tilde\psi_{k+1}, \qquad k\in\mathbb{N}.
\end{equation*}
We now describe the $\xi$-images of these elements.
For notational convenience, set 
$$
	z=\frac{t+t^{-1}}{2}.
$$

\begin{lemma}\label{lem:mat_4families}
For $k \in \mathbb{N}$, the map $\xi$ sends
\begin{align}
	A_k \quad & \longmapsto \quad \frac{(-1)^k}{2^k} \left( 1+z \right)^k \begin{pmatrix} 0 & 1 \\ 1 & 0 \end{pmatrix}, \label{eq:1st,A_k image}\\[2pt]
	B_k \quad & \longmapsto \quad \frac{(-1)^k}{2^k} \left(1+z \right)^k \begin{pmatrix} 0 & t \\ t^{-1} & 0 \end{pmatrix},  \label{eq:1st,B_k image}\\[2pt]
	\psi_{k+1} \quad & \longmapsto \quad \frac{(-1)^k}{2^{k+1}} \left(1+z \right)^k 
	\begin{pmatrix} 
	(t - t^{-1})/2 & 1 + t \\[2pt] 
	1 + t^{-1} & (t^{-1} - t)/2
	\end{pmatrix}, \label{eq:psi_k image}\\
	\tilde\psi_{k+1} \quad & \longmapsto \quad \frac{(-1)^k}{2^{k+1}} \left(1+z \right)^k 
	\begin{pmatrix} 
	(t^{-1} - t)/2 & 1 + t \\[2pt] 
	1 + t^{-1} & (t - t^{-1})/2
	\end{pmatrix}. \notag
\end{align}
\end{lemma}
\begin{proof}
By induction on $k$, using \eqref{cor:0312 psi}--\eqref{cor:0312 B} and \eqref{eq:rec tpsi} along with \eqref{eq:image of A, B}.
\end{proof}

\noindent
Next, we describe the $\xi$-images of the elements in \eqref{eq:3rd basis O}.

\begin{lemma}[{\cite[Definition~3.7]{2023TerCAlg}}]
For $k \in \mathbb{N}$, the map $\xi$ sends
\begin{align}
	W_{-k} \quad & \longmapsto \quad z^k \begin{pmatrix} 0 & 1 \\ 1 & 0 \end{pmatrix}, \label{eq:mat W_(-k)}\\[2pt]
	W_{k+1} \quad & \longmapsto \quad  z^k \begin{pmatrix} 0 & t \\ t^{-1} & 0 \end{pmatrix}, \label{eq:mat W_(k)} \\[2pt]
	\tilde{G}_{k+1} \quad & \longmapsto \quad  z^k \begin{pmatrix} t^{-1} - t & 0 \\ 0 & t - t^{-1} \end{pmatrix}. \label{eq:mat G_(k+1)}
\end{align}
\end{lemma}

\noindent
For the rest of the paper, we identify $\mathcal{O}$ with its $\xi$-image in $L$ via the injection $\xi$.

\section{Transition matrices}\label{Sec:TM}
We continue our discussion of the three bases for $\mathcal{O}$ given in \eqref{eq:1st basis O}, \eqref{eq:2nd basis O}, and \eqref{eq:3rd basis O}.
In this section, we describe the transition matrices among these three bases.
First, we consider the bases \eqref{eq:1st basis O} and \eqref{eq:2nd basis O}.
Equation \eqref{eq:tilde psi_i} gives the transition matrices between them in both directions.
Next, we consider the bases \eqref{eq:1st basis O} and \eqref{eq:3rd basis O}.

\begin{proposition}\label{lem:tm 3rd 1st}
The transition matrix from the basis \eqref{eq:3rd basis O} to the basis \eqref{eq:1st basis O} is given as follows.
For $k \in \mathbb{N}$,
\begin{align}
	A_k  & = \frac{(-1)^k}{2^k} \sum_{j=0}^k {k \choose j} W_{-j}, \label{eq(1):tm 3rd 1st} \\
	B_k  & = \frac{(-1)^k}{2^k} \sum_{j=0}^k {k \choose j} W_{j+1}, \label{eq(2):tm 3rd 1st} \\
	\psi_{k+1}  
		& = 	\frac{(-1)^k}{2^{k+1}} \sum_{j=0}^k {k \choose j} W_{-j}
			+ \frac{(-1)^k}{2^{k+1}} \sum_{j=0}^k {k \choose j} W_{j+1}
			+ \frac{(-1)^{k+1}}{2^{k+2}} \sum_{j=0}^k {k \choose j} \tilde{G}_{j+1}. \label{eq(3):tm 3rd 1st}
\end{align}
\end{proposition}

\begin{proof}
By the binomial theorem,
\begin{equation}\label{pf_eq:binom}
	(1+z)^k=\sum_{j=0}^k {k\choose j}z^j.
\end{equation}
We first show \eqref{eq(1):tm 3rd 1st}.
By \eqref{eq:1st,A_k image} and \eqref{pf_eq:binom}, we have
\begin{align*}
	A_k 
	& =  \frac{(-1)^k}{2^k} \sum_{j=0}^k {k \choose j} z^j \begin{pmatrix} 0 & 1 \\ 1 & 0 \end{pmatrix}.
\end{align*}
By this and \eqref{eq:mat W_(-k)}, we get \eqref{eq(1):tm 3rd 1st}.
Similarly, by \eqref{eq:1st,B_k image}, \eqref{eq:mat W_(k)}, and \eqref{pf_eq:binom}, we get \eqref{eq(2):tm 3rd 1st}.
Next, we show \eqref{eq(3):tm 3rd 1st}.
By \eqref{eq:psi_k image} and \eqref{pf_eq:binom}, we have
\begin{equation}\label{pf_eq(1):lem psi_k+1}
	\psi_{k+1} = \frac{(-1)^k}{2^{k+1}} \sum_{j=0}^k {k \choose j} z^j 
	\begin{pmatrix} 
	(t - t^{-1})/2 & 1 + t \\[2pt] 
	1 + t^{-1} & (t^{-1} - t)/2
	\end{pmatrix}.
\end{equation}
Observe that
\begin{equation}\label{pf_eq(3):lem psi_k+1}
	\begin{pmatrix} 
	(t - t^{-1})/2 & 1 + t \\[2pt] 
	1 + t^{-1} & (t^{-1} - t)/2
	\end{pmatrix}
	= 
	\begin{pmatrix} 
	0 & 1 \\[2pt]
	1 & 0
	\end{pmatrix}
	+
	\begin{pmatrix} 
	0 & t \\[2pt] 
	t^{-1} & 0
	\end{pmatrix}
	-
	\frac{1}{2}\begin{pmatrix} 
	t^{-1} - t & 0 \\[2pt] 
	0 & t - t^{-1}
	\end{pmatrix}.				
\end{equation}
Substitute \eqref{pf_eq(3):lem psi_k+1} into \eqref{pf_eq(1):lem psi_k+1}, and evaluate the result using \eqref{eq:mat W_(-k)}--\eqref{eq:mat G_(k+1)}.
The result follows.
\end{proof}

\begin{example}\label{ex:tm_1st_3rd, Ak}
For $0 \leq k \leq 5$, we express $A_k$  in terms of $W_0, W_{-1}, \ldots, W_{-k}$.
\begin{align*}
	A_0 & = W_0, \\[2pt]
	A_1 & = -\frac{W_0 + W_{-1}}{2}, \\[2pt]
	A_2 & = \frac{W_0 + 2 W_{-1} + W_{-2}}{4}, \\[2pt]
	A_3 & = -\frac{W_0 + 3 W_{-1} + 3 W_{-2} + W_{-3}}{8}, \\[2pt]
	A_4 & = \frac{W_0 + 4 W_{-1} + 6 W_{-2} + 4 W_{-3} + W_{-4}}{16}, \\[2pt]
	A_5 & = -\frac{W_0 + 5 W_{-1} + 10 W_{-2} + 10 W_{-3} + 5 W_{-4} + W_{-5}}{32}.
\end{align*}
\end{example}

\begin{example}\label{ex:tm_1st_3rd, Bk}
For $0 \leq k \leq 5$, we express $B_k$ in terms of $W_1, W_2, \ldots, W_{k+1}$.
\begin{align*}
	B_0 & = W_1, \\[2pt]
	B_1 & = -\frac{W_1 + W_2}{2}, \\[2pt]
	B_2 & = \frac{W_1 + 2 W_{2} + W_{3}}{4}, \\[2pt]
	B_3 & = -\frac{W_1 + 3 W_{2} + 3 W_{3} + W_{4}}{8}, \\[2pt]
	B_4 & = \frac{W_1 + 4 W_{2} + 6 W_{3} + 4 W_{4} + W_{5}}{16}, \\[2pt]
	B_5 & = -\frac{W_1 + 5 W_{2} + 10 W_{3} + 10 W_{4} + 5 W_{5} + W_{6}}{32}.
\end{align*}
\end{example}

\begin{example}\label{ex:tm_1st_3rd, psik}
For $0 \leq k \leq 4$, we express $\psi_{k+1}$ in terms of $W_0, W_{-1}, \ldots, W_{-k}$ and $W_1, W_2, \ldots, W_{k+1}$ and $\tilde G_1, \tilde G_2, \ldots, \tilde G_{k+1}$.
\begin{align*}
	\psi_1 & = \frac{W_0}{2} + \frac{W_1}{2} - \frac{\tilde G_1}{4}, \\[2pt]
	\psi_2 & = -\frac{W_0 + W_{-1}}{4} - \frac{W_1 + W_2}{4} + \frac{\tilde G_1 + \tilde G_2}{8}, \\[2pt]
	\psi_3 & = \frac{W_0 + 2W_{-1} + W_{-2}}{8} + \frac{W_1 + 2W_2 + W_3}{8} - \frac{\tilde G_1 + 2\tilde G_2 + \tilde G_3}{16}, \\[2pt]
	\psi_4 & = -\frac{W_0 + 3W_{-1} + 3W_{-2} + W_{-3}}{16} - \frac{W_1 + 3W_2 + 3W_3 + W_4}{16} + \frac{\tilde G_1 + 3\tilde G_2 + 3\tilde G_3 + \tilde G_4}{32}, \\[2pt]
	\psi_5 & = \frac{W_0 + 4W_{-1} + 6W_{-2} + 4W_{-3} + W_{-4}}{32} 
	+ \frac{W_1 + 4W_2 + 6W_3 + 4W_4 + W_5}{32} \\
	& \qquad \qquad - \frac{\tilde G_1 + 4\tilde G_2 + 6\tilde G_3 + 4\tilde G_4 + \tilde G_5}{64}.
\end{align*}
\end{example}

\begin{proposition}\label{lem:tm 1st 3rd}
The transition matrix from the basis \eqref{eq:1st basis O} to the basis \eqref{eq:3rd basis O} is given as follows.
For $k \in \mathbb{N}$,
\begin{align}
	W_{-k} & = (-1)^k \sum_{j=0}^k  {k \choose j} 2^j A_{j}, \label{eq(1):tm 1st 3rd} \\
	W_{k+1} & = (-1)^k \sum_{j=0}^k  {k \choose j} 2^j B_{j}, \label{eq(2):tm 1st 3rd} \\
	\tilde{G}_{k+1}
	& = 	(-1)^{k} \sum_{j=0}^k  {k \choose j} 2^{j+1} A_j
		+ (-1)^{k} \sum_{j=0}^k  {k \choose j} 2^{j+1} B_j
		+ (-1)^{k+1} \sum_{j=0}^k  {k \choose j} 2^{j+2} \psi_{j+1}. \label{eq(3):tm 1st 3rd}
\end{align}
\end{proposition}

\begin{proof}
By the binomial theorem,
\begin{equation}\label{pf_eq:inverse binom}
	z^k=((1+z)-1)^k
	=\sum_{j=0}^k {k\choose j}(-1)^{k-j}(1+z)^j.
\end{equation}
We first show \eqref{eq(1):tm 1st 3rd}.
By \eqref{eq:mat W_(-k)} and \eqref{pf_eq:inverse binom},
\begin{align*}
	W_{-k}
	&=\sum_{j=0}^k {k\choose j}(-1)^{k-j}(1+z)^j
	\begin{pmatrix}
	0&1\\
	1&0
	\end{pmatrix}
	=(-1)^k\sum_{j=0}^k {k\choose j}2^j
	\frac{(-1)^j}{2^j}(1+z)^j
	\begin{pmatrix}
	0&1\\
	1&0
	\end{pmatrix}.
\end{align*}
Hence, by \eqref{eq:1st,A_k image}, we get \eqref{eq(1):tm 1st 3rd}.
Similarly, by \eqref{eq:mat W_(k)}, \eqref{pf_eq:inverse binom}, and \eqref{eq:1st,B_k image}, we get \eqref{eq(2):tm 1st 3rd}.
Next, we show \eqref{eq(3):tm 1st 3rd}.
By \eqref{eq:mat G_(k+1)} and \eqref{pf_eq:inverse binom},
\begin{equation}\label{pf_eq(1):lem 1st 3rd}
	\tilde G_{k+1}
	=
	\sum_{j=0}^k {k\choose j}(-1)^{k-j}(1+z)^j
	\begin{pmatrix}
	t^{-1}-t&0\\
	0&t-t^{-1}
	\end{pmatrix}.
\end{equation}
Observe that
\begin{equation}\label{pf_eq(2):lem 1st 3rd}
	\begin{pmatrix}
	t^{-1}-t&0\\[2pt]
	0&t-t^{-1}
	\end{pmatrix}
	=
	2\begin{pmatrix}
	0&1\\[2pt]
	1&0
	\end{pmatrix}
	+
	2\begin{pmatrix}
	0&t\\[2pt]
	t^{-1}&0
	\end{pmatrix}
	-
	2\begin{pmatrix}
	(t-t^{-1})/2&1+t\\[2pt]
	1+t^{-1}&(t^{-1}-t)/2
	\end{pmatrix}.
\end{equation}
Substitute \eqref{pf_eq(2):lem 1st 3rd} into \eqref{pf_eq(1):lem 1st 3rd}, and evaluate the result using \eqref{eq:1st,A_k image}--\eqref{eq:psi_k image}.
The result follows.
\end{proof}

\begin{example}\label{ex:tm_3rd_1st, Ak}
For $0 \leq k \leq 5$, we express $W_{-k}$  in terms of $A_0, A_1, \ldots, A_{k}$.
\begin{align*}
	W_0 & = A_0, \\[2pt]
	W_{-1} & = - (A_0 + 1 \cdot 2 A_1 ), \\[2pt]
	W_{-2} & =  A_0 + 2\cdot 2 A_1 + 1 \cdot 2^2 A_2, \\[2pt]
	W_{-3} & = -  ( A_0 + 3 \cdot 2 A_1 + 3 \cdot 2^2 A_2 + 1\cdot 2^3 A_3), \\[2pt]
	W_{-4} & = A_0 + 4 \cdot 2 A_1 + 6 \cdot 2^2 A_2 + 4 \cdot 2^3 A_3 + 1 \cdot 2^4 A_4, \\[2pt]
	W_{-5} & = -(A_0 + 5 \cdot 2 A_1 + 10 \cdot 2^2 A_2 + 10 \cdot 2^3 A_3 + 5 \cdot 2^4 A_4 + 1 \cdot 2^5 A_5).
\end{align*}
\end{example}

\begin{example}\label{ex:tm_3rd_1st, Bk}
For $0 \leq k \leq 5$, we express $W_{k+1}$ in terms of $B_0, B_1, \ldots, B_{k}$.
\begin{align*}
	W_1 & = B_0, \\[2pt]
	W_{2} & = - (B_0 + 1 \cdot 2 B_1 ), \\[2pt]
	W_{3} & =  B_0 + 2\cdot 2 B_1 + 1 \cdot 2^2 B_2, \\[2pt]
	W_{4} & = -  ( B_0 + 3 \cdot 2 B_1 + 3 \cdot 2^2 B_2 + 1\cdot 2^3 B_3), \\[2pt]
	W_{5} & = B_0 + 4 \cdot 2 B_1 + 6 \cdot 2^2 B_2 + 4 \cdot 2^3 B_3 + 1 \cdot 2^4 B_4, \\[2pt]
	W_{6} & = -(B_0 + 5 \cdot 2 B_1 + 10 \cdot 2^2 B_2 + 10 \cdot 2^3 B_3 + 5 \cdot 2^4 B_4 + 1 \cdot 2^5 B_5).
\end{align*}
\end{example}

\begin{example}\label{ex:tm_3rd_1st, psi_k}
For $0 \leq k \leq 4$, we express $\tilde G_{k+1}$ in terms of $A_0, A_1, \ldots, A_k$ and $B_0, B_1, \ldots, B_{k}$ and $\psi_1, \psi_2, \ldots, \psi_{k+1}$.
\begin{align*}
	\tilde G_1 & = 2A_0 + 2B_0 - 2^2 \psi_1, \\[2pt]
	\tilde G_2 & = -1\cdot 2A_0 - 1 \cdot 2^2 A_1 - 1\cdot 2 B_0 - 1\cdot 2^2 B_1 + 1\cdot 2^2 \psi_1 + 1\cdot 2^3 \psi_2, \\[2pt]
	\tilde G_3 & = 1\cdot 2 A_0 + 2\cdot2^2 A_1 + 1\cdot 2^3 A_2 
			+ 1 \cdot 2 B_0 + 2 \cdot 2^2 B_1 + 1 \cdot 2^3 B_2 \\[2pt]
	& \qquad \qquad - 1 \cdot 2^2 \psi_1 - 2 \cdot 2^3 \psi_2 - 1 \cdot 2^4 \psi_3, \\[2pt]
	\tilde G_4 & = -1 \cdot 2 A_0 - 3 \cdot 2^2 A_1 - 3 \cdot 2^3 A_2 - 1 \cdot 2^4 A_3
			- 1 \cdot 2 B_0 - 3 \cdot 2^2 B_1 - 3 \cdot 2^3 B_2 - 1\cdot 2^4 B_3 \\[2pt]
	& \qquad \qquad	
			+ 1 \cdot 2^2 \psi_1 + 3 \cdot 2^3 \psi_2 + 3 \cdot 2^4 \psi_3 + 1\cdot 2^5 \psi_4, \\[2pt]
	\tilde G_5 & = 1 \cdot 2 A_0 + 4 \cdot 2^2 A_1 + 6 \cdot 2^3 A_2 + 4 \cdot 2^4 A_3 + 1 \cdot 2^5 A_4 \\[2pt]
			& \qquad \qquad + 1 \cdot 2 B_0 + 4 \cdot 2^2 B_1 + 6 \cdot 2^3 B_2 + 4 \cdot 2^4 B_3 + 1 \cdot 2^5 B_4\\[2pt]
	& \qquad \qquad \qquad \qquad
			- 1 \cdot 2^2 \psi_1 - 4 \cdot 2^3 \psi_2 - 6 \cdot 2^4 \psi_3 - 4 \cdot 2^5 \psi_4 - 1 \cdot 2^6 \psi_5.
\end{align*}
\end{example}

Next, we consider the transition matrices between the bases \eqref{eq:2nd basis O} and \eqref{eq:3rd basis O}.

\begin{proposition}\label{lem:tm 3rd 2nd}
The transition matrix from the basis \eqref{eq:3rd basis O} to the basis \eqref{eq:2nd basis O} is given as follows.
For $k \in \mathbb{N}$,
\begin{align*}
	A_k  & = \frac{(-1)^k}{2^k} \sum_{j=0}^k {k \choose j} W_{-j}, \\
	B_k  & = \frac{(-1)^k}{2^k} \sum_{j=0}^k {k \choose j} W_{j+1}, \\
	\tilde\psi_{k+1}  
		& = 	\frac{(-1)^k}{2^{k+1}} \sum_{j=0}^k {k \choose j} W_{-j}
			+ \frac{(-1)^k}{2^{k+1}} \sum_{j=0}^k {k \choose j} W_{j+1}
			+ \frac{(-1)^k}{2^{k+2}} \sum_{j=0}^k {k \choose j} \tilde{G}_{j+1}.
\end{align*}
\end{proposition}
\begin{proof}
The result follows from Proposition \ref{lem:tm 3rd 1st} and \eqref{eq:tilde psi_i}.
\end{proof}

\begin{example}
For $0 \leq k \leq 4$, we express $\tilde\psi_{k+1}$ in terms of $W_0, W_{-1}, \ldots, W_{-k}$ and $W_1, W_2, \ldots, W_{k+1}$ and $\tilde G_1, \tilde G_2, \ldots, \tilde G_{k+1}$.
\begin{align*}
	\tilde\psi_1 & = \frac{W_0}{2} + \frac{W_1}{2} + \frac{\tilde G_1}{4}, \\[2pt]
	\tilde\psi_2 & = -\frac{W_0 + W_{-1}}{4} - \frac{W_1 + W_2}{4} - \frac{\tilde G_1 + \tilde G_2}{8}, \\[2pt]
	\tilde\psi_3 & = \frac{W_0 + 2W_{-1} + W_{-2}}{8} + \frac{W_1 + 2W_2 + W_3}{8} + \frac{\tilde G_1 + 2\tilde G_2 + \tilde G_3}{16}, \\[2pt]
	\tilde\psi_4 & = -\frac{W_0 + 3W_{-1} + 3W_{-2} + W_{-3}}{16} - \frac{W_1 + 3W_2 + 3W_3 + W_4}{16}  - \frac{\tilde G_1 + 3\tilde G_2 + 3\tilde G_3 + \tilde G_4}{32}, \\[2pt]
	\tilde\psi_5 & = \frac{W_0 + 4W_{-1} + 6W_{-2} + 4W_{-3} + W_{-4}}{32} 
	+ \frac{W_1 + 4W_2 + 6W_3 + 4W_4 + W_5}{32} \\
	& \qquad \qquad + \frac{\tilde G_1 + 4\tilde G_2 + 6\tilde G_3 + 4\tilde G_4 + \tilde G_5}{64}.
\end{align*}
\end{example}

\begin{proposition}\label{lem:tm 2nd 3rd}
The transition matrix from the basis \eqref{eq:2nd basis O} to the basis \eqref{eq:3rd basis O} is given as follows.
For $k \in \mathbb{N}$,
\begin{align*}
	W_{-k} & = (-1)^k \sum_{j=0}^k  {k \choose j} 2^j A_{j}, \\
	W_{k+1} & = (-1)^k \sum_{j=0}^k  {k \choose j} 2^j B_{j}, \\
	\tilde{G}_{k+1}
	& = 	(-1)^{k+1} \sum_{j=0}^k  {k \choose j} 2^{j+1} A_j
		+ (-1)^{k+1} \sum_{j=0}^k  {k \choose j} 2^{j+1} B_j
		+ (-1)^{k} \sum_{j=0}^k  {k \choose j} 2^{j+2} \tilde\psi_{j+1}. 
\end{align*}
\end{proposition}
\begin{proof}
The result follows from Proposition \ref{lem:tm 1st 3rd} and \eqref{eq:tilde psi_i}.
\end{proof}

\begin{example}
For $0 \leq k \leq 4$, we express $\tilde G_{k+1}$ in terms of $A_0, A_1, \ldots, A_k$ and $B_0, B_1, \ldots, B_{k}$ and $\tilde\psi_1, \tilde\psi_2, \ldots, \tilde\psi_{k+1}$.
\begin{align*}
	\tilde G_1 & = - 2A_0 - 2B_0 + 2^2 \tilde\psi_1, \\[2pt]
	\tilde G_2 & = 1\cdot 2A_0 + 1 \cdot 2^2 A_1 + 1\cdot 2 B_0 + 1\cdot 2^2 B_1 - 1\cdot 2^2 \tilde\psi_1 - 1\cdot 2^3 \tilde\psi_2, \\[2pt]
	\tilde G_3 & = -1\cdot 2 A_0 - 2\cdot2^2 A_1 - 1\cdot 2^3 A_2 
			- 1 \cdot 2 B_0 - 2 \cdot 2^2 B_1 - 1\cdot 2^3 B_2 \\[2pt]
	& \qquad \qquad + 1 \cdot 2^2 \tilde\psi_1 + 2 \cdot 2^3 \tilde\psi_2 + 1 \cdot 2^4\tilde\psi_3, \\[2pt]
	\tilde G_4 & = 1 \cdot 2 A_0 + 3 \cdot 2^2 A_1 + 3 \cdot 2^3 A_2 + 1\cdot 2^4 A_3
			+ 1 \cdot 2 B_0 + 3 \cdot 2^2 B_1 + 3 \cdot 2^3 B_2 + 1\cdot 2^4 B_3 \\[2pt]
	& \qquad \qquad	
			- 1 \cdot 2^2 \tilde\psi_1 - 3 \cdot 2^3 \tilde\psi_2 - 3 \cdot 2^4 \tilde\psi_3 - 1\cdot 2^5 \tilde\psi_4, \\[2pt]
	\tilde G_5 & = -1 \cdot 2 A_0 - 4 \cdot 2^2 A_1 - 6 \cdot 2^3 A_2 - 4 \cdot 2^4 A_3 - 1 \cdot 2^5 A_4 \\[2pt]
			& \qquad \qquad -1 \cdot 2 B_0 - 4 \cdot 2^2 B_1 - 6 \cdot 2^3 B_2 - 4 \cdot 2^4 B_3 - 1 \cdot 2^5 B_4\\[2pt]
	& \qquad \qquad \qquad \qquad
			+ 1 \cdot 2^2 \tilde\psi_1 + 4 \cdot 2^3 \tilde\psi_2 + 6 \cdot 2^4 \tilde\psi_3 + 4 \cdot 2^5 \tilde\psi_4 + 1 \cdot 2^6 \tilde\psi_5.
\end{align*}

\end{example}

\section{Tridiagonal pairs and tridiagonal systems}\label{Sec:TD_pair}

For the rest of this paper, we consider the finite-dimensional irreducible $\mathcal{O}$-modules.
Our goal is to describe how the following elements of $\mathcal{O}$ act on these modules:
\begin{equation}\label{eq:module elements}
	A_k, \qquad B_k, \qquad \psi_{k+1}, \qquad \tilde\psi_{k+1}, \qquad k\in\mathbb{N}.
\end{equation}
To this end, we use the theory of tridiagonal pairs.
We begin by recalling the definition of a tridiagonal pair.

\begin{definition}[{\cite[Definition~1.1]{ITT2001}}]\label{def1}
Let $V$ denote a vector space of finite positive dimension. 
By a \emph{tridiagonal pair} (or \emph{TD pair}) on $V$, we mean an ordered pair $A, A^*$ of elements of $\operatorname{End}(V)$ that satisfy (i)--(iv) below.
\begin{enumerate}[label=(\roman*),font=\rm,itemsep=0mm]
	\item Each of $A$ and $A^*$ is diagonalizable on $V$. 
	\item There exists an ordering $\{ V_i \}_{i=0}^d$ of the eigenspaces of $A$ such that 
	\begin{equation}\label{eq1:TDpairDef}
	A^*V_i \subseteq V_{i-1}+V_i+V_{i+1} \qquad (0 \leq i \leq d),
	\end{equation}
	where $V_{-1}=V_{d+1}=0$. 
	\item There exists an ordering $\{ V^*_i \}_{i=0}^{\delta}$ of the eigenspaces of $A^*$ such that 
	\begin{equation}\label{eq2:TDpairDef}
	AV_i^* \subseteq V^*_{i-1}+V^*_i+V^*_{i+1} \qquad (0 \leq i \leq \delta),
	\end{equation}
	where $V^*_{-1}=V^*_{\delta+1}=0$. 
	\item There does not exist a subspace $W$ of $V$ such that $W \ne 0$, $W \ne V$, $A W \subseteq W$, $A^* W \subseteq W$. 
\end{enumerate}
\end{definition} 

\begin{note}\label{note:A*}
It is common notation to use $A^*$ to denote the conjugate transpose of $A$.  
We are not using this notation in this paper.  
In a tridiagonal pair $A, A^*$, the $\mathbb{F}$-linear maps $A$ and $A^*$ are arbitrary subject to conditions (i)–(iv) in Definition~\ref{def1}.
\end{note}

Let  $A, A^*$ denote a TD pair on $V$ as in Definition~\ref{def1}.  
By \cite[Lemma 4.5]{ITT2001} the integers $d$ and $\delta$ from \eqref{eq1:TDpairDef} and \eqref{eq2:TDpairDef} are equal. 
We call this common value the \emph{diameter} of the TD pair. 
An ordering $\{V_i\}_{i=0}^d$ of the eigenspaces of $A$ is called \emph{standard} whenever it satisfies \eqref{eq1:TDpairDef}.
Let $\{V_i\}_{i=0}^d$ be a standard ordering of the eigenspaces of $A$.
Then $\{V_{d-i}\}_{i=0}^d$ is also standard, and no further ordering is standard \cite[Lemma 2.4]{ITT2001}.
A similar result holds for the eigenspaces of $A^*$.
We now consider the primitive idempotents of $A$ and $A^*$.
An ordering of the primitive idempotents of $A$ (resp. $A^*$) is called \emph{standard} whenever the corresponding ordering of the eigenspaces of $A$ (resp. $A^*$) is standard.

\smallskip
When working with a TD pair, it is convenient to use a closely related concept known as a TD system.

\begin{definition}[{\cite[Definition~2.1]{ITT2001}}]\label{def:TDsystem}
Let $V$ denote a vector space of finite positive dimension. 
By a \emph{tridiagonal system} (or \emph{TD system}) on $V$, we mean a sequence 
\begin{equation*}
	\Phi = ( A; \{ E_i \}_{i=0}^d; A^*; \{ E_i^* \}_{i=0}^{d} )
\end{equation*}
that satisfies (i)--(iii) below.
\begin{enumerate}[label=(\roman*),font=\rm,itemsep=0mm]
	\item The pair $A$, $A^*$ is a TD pair on $V$. 
	\item $\{E_i\}_{i=0}^d$ is a standard ordering of the primitive idempotents of $A$. 
	\item $\{E^*_i\}_{i=0}^d$ is a standard ordering of the primitive idempotents of $A^*$.
\end{enumerate}
\end{definition}

Let $\Phi=(A;\{E_i\}_{i=0}^d;A^*;\{E_i^*\}_{i=0}^d)$ denote a TD system on $V$.
For $0\leq i\leq d$, let $\theta_i$ denote the eigenvalue of $A$ associated with $E_i$, and let $\theta_i^*$ denote the eigenvalue of $A^*$ associated with $E_i^*$.
We call $\{\theta_i\}_{i=0}^d$ the \emph{eigenvalue sequence} of $\Phi$ and $\{\theta_i^*\}_{i=0}^d$ the \emph{dual eigenvalue sequence} of $\Phi$.

\smallskip
Let $A, A^*$ denote a TD pair on $V$.
Let $\{E_i\}_{i=0}^d$ and $\{E_i^*\}_{i=0}^d$ be standard orderings of the primitive idempotents of $A$ and $A^*$, respectively.
By construction, each of the following is a TD system on $V$:
\begin{align*}
	& ( A; \{ E_i \}_{i=0}^d; A^*; \{ E_i^* \}_{i=0}^{d} ), \\
	& ( A; \{ E_i \}_{i=0}^d; A^*; \{ E_{d-i}^* \}_{i=0}^{d} ), \\
	& ( A; \{ E_{d-i} \}_{i=0}^d; A^*; \{ E_i^* \}_{i=0}^{d} ), \\
	& ( A; \{ E_{d-i} \}_{i=0}^d; A^*; \{ E_{d-i}^* \}_{i=0}^{d} ).
\end{align*}
By an eigenvalue sequence (resp. dual eigenvalue sequence) of $A, A^*$, we mean an eigenvalue sequence (resp. dual eigenvalue sequence) of one of the TD systems listed above.
For more information on TD pairs and TD systems, see \cite{ITT2001, INP2011LAA, IT2004JPAA, IT2007RJ, IT2007LAA, NT2007LAA, NT2008LAA2, NT2008LAA}.

\section{Decompositions and flags}\label{Sec:DecompFlag}

We continue to discuss the TD system $\Phi=( A; \{ E_i \}_{i=0}^d; A^*; \{ E_i^* \}_{i=0}^{d} )$ on $V$.
In this section, we discuss some decompositions and flags associated with $\Phi$.
We begin by recalling the notion of a decomposition.
By a \emph{decomposition} of $V$, we mean a sequence $\{U_i\}_{i=0}^d$ of nonzero subspaces whose direct sum is $V$.
For example, each of the standard orderings of the eigenspaces of $A$ and $A^*$ is a decomposition of $V$.
Let $\{U_i\}_{i=0}^d$ denote a decomposition of $V$.
For $0 \leq i \leq d$ we call $U_i$ the \emph{$i$th component} of the decomposition.
By construction, if $\{U_i\}_{i=0}^d$ is a decomposition of $V$, then so is $\{U_{d-i}\}_{i=0}^d$.

\begin{example}[{\cite[Section~4]{ITT2001}}]
Let $\Phi=( A; \{ E_i \}_{i=0}^d; A^*; \{ E_i^* \}_{i=0}^{d} )$ be a TD system on $V$.
For $0 \leq i \leq d$ define
\begin{equation}\label{eq:split U_i}
	U_i = (E_0^*V + E_1^*V + \cdots + E_i^*V) \cap (E_0V + E_1V + \cdots + E_{d-i}V).
\end{equation}
Note that $U_0 = E_0^*V$ and $U_d = E_0V$.
Then $\{U_i\}_{i=0}^d$ is a decomposition of $V$.
This decomposition is called a \emph{split} decomposition of $V$.
\end{example}

\noindent
We recall some facts about the split decomposition $\{U_i\}_{i=0}^d$.
By \cite[Theorem~4.6]{ITT2001}, for $0 \leq i \leq d$ both
\begin{align*}
	U_0 + U_1 + \cdots + U_i & = E_0^*V + E_1^*V + \cdots + E_i^*V, \\
	U_i + U_{i+1} + \cdots + U_d & = E_0V + E_{1}V + \cdots + E_{d-i}V.
\end{align*}
Moreover, 
\begin{align*}
	&& (A - \theta_{d-i} I )U_i & \subseteq U_{i+1} \qquad (0 \leq i \leq d-1), && (A - \theta_0I)U_d = 0, && \\
	&& (A^* - \theta_i^* I )U_i & \subseteq U_{i-1} \qquad (1 \leq i \leq d), && (A^* - \theta_0^*I)U_0 = 0. &&
\end{align*}

\noindent
We now describe six decompositions of $V$ associated with $\Phi$.
These decompositions will be used throughout the rest of the paper.
In each row of the table below, we list one of these decompositions; see \cite[Lemma~4.2]{IT2007RJ}. 
\begin{equation}\label{table:6decomp}
\renewcommand{\arraystretch}{1.2} 
\begin{tabular}{c|c}
	decomposition name & $i$th component of the decomposition \\
	\hline \hline
	$[0, d]$ & $E_iV$ \\
	$[0^*, d^*]$ & $E_i^*V$ \\
	$[0^*, 0]$ & $(E_0^*V + \cdots + E_i^*V) \cap (E_0V + \cdots + E_{d-i}V)$ \\
	$[0^*, d]$ & $(E_0^*V + \cdots + E_i^*V) \cap (E_iV + \cdots + E_dV)$ \\
	$[d^*, 0]$ & $(E_{d-i}^*V + \cdots + E_d^*V) \cap (E_0V + \cdots + E_{d-i}V)$ \\
	$[d^*, d]$ & $(E_{d-i}^*V + \cdots + E_d^*V) \cap (E_iV + \cdots + E_dV)$ \\
\end{tabular}
\end{equation}
Observe that the decomposition $[0^*,0]$ in \eqref{table:6decomp} is the split decomposition $\{U_i\}_{i=0}^d$ from \eqref{eq:split U_i}.

\smallskip
Next, we recall the notion of a flag.
By a \emph{flag on $V$}, we mean a sequence $\{F_i\}_{i=0}^d$ of subspaces of $V$ such that
$$
	F_{i-1}\subsetneq F_i \qquad (0\leq i\leq d),
$$
where $F_{-1}=0$ and $F_d=V$.
For $0 \leq i \leq d$ we call $F_i$ the \emph{$i$th component} of the flag.
A decomposition of $V$ naturally yields a flag on $V$ as follows.
Let $\{U_i\}_{i=0}^d$ denote a decomposition of $V$.
For $0 \leq i \leq d$, define $F_i = U_0 + U_1 + \cdots + U_i$.
Then the sequence $\{F_i\}_{i=0}^d$ is a flag on $V$; we say that $\{F_i\}_{i=0}^d$ is \emph{induced} by the decomposition $\{U_i\}_{i=0}^d$.
Let $\{F_i\}_{i=0}^d$ and $\{F'_i\}_{i=0}^d$ denote flags on $V$.
These flags are said to be \emph{opposite} whenever there exists a decomposition $\{U_i\}_{i=0}^d$ of $V$ such that $\{U_i\}_{i=0}^d$ induces $\{F_i\}_{i=0}^d$ and $\{U_{d-i}\}_{i=0}^d$ induces $\{F'_i\}_{i=0}^d$. 
In this case, we have $U_i = F_i \cap F'_{d-i}$ for $0 \leq i \leq d$.
The TD system $\Phi$ gives rise to four flags on $V$.
The table below displays these four flags; see \cite[Definition~4.2]{2023TerJCTA}. 
\begin{equation}\label{table:4flags}
\renewcommand{\arraystretch}{1.2} 
\begin{tabular}{c|c}
	flag name & $i$th component of the flag \\
	\hline \hline
	$[0]$ & $E_0V + E_1V + \cdots + E_iV$ \\
	$[d]$ & $E_dV + E_{d-1}V + \cdots + E_{d-i}V$ \\
	$[0^*]$ & $E_0^*V + E_1^*V + \cdots + E_i^*V$ \\
	$[d^*]$ & $E_d^*V + E_{d-1}^*V + \cdots + E_{d-i}^*V$
\end{tabular}
\end{equation}
By \cite[Lemma~4.3]{2023TerJCTA}, the four flags in \eqref{table:4flags} are mutually opposite.

\smallskip
We finish this section with some terminology and a lemma concerning decompositions and flags.
Let $\{U_i\}_{i=0}^d$ denote a decomposition of $V$.
Let $\{F_i\}_{i=0}^d$ denote the flag on $V$ induced by $\{U_i\}_{i=0}^d$ and let $\{F'_i\}_{i=0}^d$ denote the flag on $V$ induced by $\{U_{d-i}\}_{i=0}^d$; that is,
$$
	F_i = U_0 + \cdots + U_i, \qquad \qquad
	F'_i = U_{d-i} + \cdots + U_d \qquad \qquad
	(0 \leq i \leq d).
$$
Let $\rho \in \operatorname{End}(V)$.
We say that $\rho$ \emph{stabilizes} $\{F_i\}_{i=0}^d$ whenever $\rho F_i \subseteq F_i$ for $0 \leq i \leq d$.
We say that $\rho$ \emph{raises} $\{F_i\}_{i=0}^d$ whenever $\rho F_i \subseteq F_{i+1}$ for $0 \leq i \leq d-1$.
For notational convenience, set $U_{-1}=U_{d+1}=0$.

\begin{lemma}[{\cite[Lemmas~6.3--6.6]{2023TerJCTA}}]\label{lem:raise,stabilizes}
With the above notation, let $\rho \in \operatorname{End}(V)$. 
Then the following {\rm(i)}--{\rm(iv)} hold.
\begin{itemize}
\item[\rm(i)] $\rho$ raises $\{F_i\}_{i=0}^d$ and $\{F'_i\}_{i=0}^d$ 
	\ if and only if \ $\rho U_i \subseteq U_{i-1} + U_i + U_{i+1}$ for $0 \leq i \leq d$.
\item[\rm(ii)] $\rho$ raises $\{F_i\}_{i=0}^d$ and stabilizes $\{F'_i\}_{i=0}^d$ 
	\ if and only if \  $\rho U_i \subseteq U_i + U_{i+1}$ for $0 \leq i \leq d$.
\item[\rm(iii)] $\rho$ stabilizes $\{F_i\}_{i=0}^d$ and raises $\{F'_i\}_{i=0}^d$ 
	\ if and only if \  $\rho U_i \subseteq U_{i-1} + U_{i}$ for $0 \leq i \leq d$.
\item[\rm(iv)] $\rho$ stabilizes $\{F_i\}_{i=0}^d$ and $\{F'_i\}_{i=0}^d$ 
	\ if and only if \  $\rho U_i \subseteq U_{i}$ for $0 \leq i \leq d$.
\end{itemize}
\end{lemma}

\section{The tetrahedron diagram}\label{Sec:decomp-diag}
In the previous section, we discussed six decompositions of $V$ associated with the TD system $\Phi$.
In this section, we illustrate these decompositions using a diagram.
This diagram will be used in Section~\ref{Sec:action-on-V} to describe how the elements in \eqref{eq:module elements} act on the six decompositions.

\smallskip
Let $\{U_i\}_{i=0}^d$ denote a decomposition of $V$.
We illustrate this decomposition using the following diagram:
\begin{center}
\begin{tikzpicture}
    \coordinate (U0) at (0,0);
    \coordinate (U1) at (1,0);
    \coordinate (U2) at (2,0);
    \coordinate (U3) at (3,0);
    \coordinate (U4) at (4,0);
    \coordinate (Ud1) at (5,0);
    \coordinate (Ud) at (6,0);

    \draw (U0) -- (U3);
    \draw[dashed] (U3) -- (U4);
    \draw (U4) -- (Ud);

    \foreach \x in {U0,U1,U2,U3,U4,Ud1,Ud}
        \fill (\x) circle (2.5pt);

    \node[below=6pt] at (U0) {$U_0$};
    \node[below=6pt] at (U1) {$U_1$};
    \node[below=6pt] at (U2) {$U_2$};
    \node[below=6pt] at (3.5,0) {$\cdots$};
    \node[below=6pt] at (Ud1) {$U_{d-1}$};
    \node[below=6pt] at (Ud) {$U_d$};
\end{tikzpicture}
\end{center}
We might omit the labels $U_i$ when they are clear from the context.

\smallskip
Let $\{U_i\}_{i=0}^d$ and $\{U_i'\}_{i=0}^d$ denote decompositions of $V$ that induce the same flag.
In other words,
$$
	U_0+\cdots+U_i=U_0'+\cdots+U_i'
	\qquad (0\leq i\leq d).
$$
We illustrate this relation using the following diagram:
\begin{center}
\begin{tikzpicture}
    \coordinate (O) at (0,0);

    \coordinate (U0) at (0,0);
    \coordinate (U1) at (1,0);
    \coordinate (U2) at (2,0);
    \coordinate (U3) at (3,0);
    \coordinate (U4) at (4,0);
    \coordinate (Ud1) at (5,0);
    \coordinate (Ud) at (6,0);

    \coordinate (Up0) at (0,0);
    \coordinate (Up1) at ({cos(30)}, {sin(30)});
    \coordinate (Up2) at ({2*cos(30)}, {2*sin(30)});
    \coordinate (Up3) at ({3*cos(30)}, {3*sin(30)});
    \coordinate (Up4) at ({4*cos(30)}, {4*sin(30)});
    \coordinate (Upd1) at ({5*cos(30)}, {5*sin(30)});
    \coordinate (Upd) at ({6*cos(30)}, {6*sin(30)});

    \draw (U0) -- (U3);
    \draw (U4) -- (Ud);
    \draw[dashed] (U3) -- (U4);
    
    \draw (Up0) -- (Up3);
    \draw (Up4) -- (Upd);
    \draw[dashed] (Up3) -- (Up4);

    \foreach \x in {U0,U1,U2,U3,U4,Ud1,Ud}
        \fill (\x) circle (2.5pt);

    \foreach \x in {Up1,Up2,Up3, Up4, Upd1,Upd}
        \fill (\x) circle (2.5pt);

    \node[below=6pt] at (U0) {$U_0$};
    \node[below=6pt] at (U1) {$U_1$};
    \node[below=6pt] at (U2) {$U_2$};
    \node[below=6pt] at (Ud1) {$U_{d-1}$};
    \node[below=6pt] at (Ud) {$U_d$};

    \node[above=6pt] at (Up0) {$U_0'$};
    \node[above=6pt] at (Up1) {$U_1'$};
    \node[above=6pt] at (Up2) {$U_2'$};
    \node[above=6pt] at (Upd1) {$U_{d-1}'$};
    \node[above=6pt] at (Upd) {$U_d'$};
    
\end{tikzpicture}
\end{center}

Recall the TD system $\Phi$ on $V$.
Associated with $\Phi$ are the six decompositions given in \eqref{table:6decomp}.
These decompositions are described in the following diagram.

\begin{center}
\scalemath{0.7}{
\begin{tikzpicture}
  [scale=1,thick,auto=left, every node/.style={circle, inner sep=2pt}] 
  \node[fill=black,label={[below=0pt]:{\Large$E_0V$}}] (b0) at (0,0) {};
  \node[fill=black,label=below:{}] (b1) at (1,0)  {};
  \node[fill=black,label=below:{}] (b2) at (1.9,0)  {};
  \node[fill=black,label=below:{}] (b3) at (2.8,0) {};
  \node[fill=black,label=below:{}] (b4) at (4.2,0)  {};
  \node[fill=black,label=below:{}] (b5) at (5.1,0)  {};
  \node[fill=black,label=below:{}] (b6) at (6,0)  {};
  \node[fill=black,label={[below=0pt]:{\Large$E_dV$}}] (b7) at (7,0)  {};
  
  \node[fill=black,label={[above=-5pt]:{\Large$E_0^*V$}}] (t0) at (0,7)  {};
  \node[fill=black,label=above:{}] (t1) at (1,7)  {};
  \node[fill=black,label=above:{}] (t2) at (1.9,7)  {};
  \node[fill=black,label=above:{}] (t3) at (2.8,7)  {};
  \node[fill=black,label=above:{}] (t4) at (4.2,7)  {};
  \node[fill=black,label=above:{}] (t5) at (5.1,7)  {}; 
  \node[fill=black,label=above:{}] (t6) at (6,7)  {}; 
  \node[fill=black,label={[above=-5pt]:{\Large$E_d^*V$}}] (t7) at (7,7)  {}; 

  \node[fill=black,label=left:{}] (l1) at (0,1)  {};
  \node[fill=black,label=left:{}] (l2) at (0,1.9)  {};
  \node[fill=black,label=left:{}] (l3) at (0,2.8)  {};
  \node[fill=black,label=left:{}] (l4) at (0,4.2)  {};
  \node[fill=black,label=left:{}] (l5) at (0,5.1)  {};
  \node[fill=black,label=left:{}] (l6) at (0,6)  {};

  \node[fill=black,label=right:{}] (r1) at (7,1)  {};
  \node[fill=black,label=right:{}] (r2) at (7,1.9)  {};
  \node[fill=black,label=right:{}] (r3) at (7,2.8)  {};
  \node[fill=black,label=right:{}] (r4) at (7,4.2)  {};
  \node[fill=black,label=right:{}] (r5) at (7,5.1)  {};
  \node[fill=black,label=right:{}] (r6) at (7,6)  {};

  \node[fill=black,label=right:{}] (d1) at (1,1)  {};
  \node[fill=black,label=right:{}] (d2) at (1.9,1.9)  {};
  \node[fill=black,label=right:{}] (d3) at (2.8,2.8)  {};
  \node[fill=black,label=right:{}] (d4) at (4.2,4.2)  {};
  \node[fill=black,label=right:{}] (d5) at (5.1,5.1)  {};
  \node[fill=black,label=right:{}] (d6) at (6,6)  {};

  \node[fill=black,label=right:{}] (e1) at (1,6)  {};
  \node[fill=black,label=right:{}] (e2) at (1.9,5.1)  {};
  \node[fill=black,label=right:{}] (e3) at (2.8,4.2)  {};
  \node[fill=black,label=right:{}] (e4) at (4.2,2.8)  {};
  \node[fill=black,label=right:{}] (e5) at (5.1,1.9)  {};
  \node[fill=black,label=right:{}] (e6) at (6,1)  {};

  \foreach \from/\to in
  {b0/b1, b1/b2, b2/b3, b4/b5, b5/b6, b6/b7,
  t0/t1, t1/t2, t2/t3, t4/t5, t5/t6, t6/t7,                  
  b0/l1, l1/l2, l2/l3,  l4/l5, l5/l6, l6/t0,
  b7/r1, r1/r2, r2/r3, r4/r5, r5/r6, r6/t7,
  t0/e1, e1/e2, e2/e3, e4/e5, e5/e6, e6/b7,
  b0/d1, d1/d2, d2/d3, d4/d5, d5/d6, d6/t7}       
      \draw (\from) -- (\to);
      \draw[dashed] (b3) -- (b4);
      \draw[dashed] (l3) -- (l4);
      \draw[dashed] (r3) -- (r4);
      \draw[dashed] (t3) -- (t4);
      \draw[dashed] (e3) -- (e4);
      \draw[dashed] (d3) -- (d4);
\end{tikzpicture}}
\end{center}

The above diagram is called the \emph{tetrahedron diagram} of $\Phi$.
The six edges of the diagram represent the six decompositions of $V$ listed in \eqref{table:6decomp}.
For more information on the tetrahedron diagram, see \cite[Section~5]{2023TerJCTA}.

\smallskip
Next, we use the tetrahedron diagram to illustrate the action of $A$ and $A^*$ on the six decompositions of $V$ in \eqref{table:6decomp}.
Let $\{U_i\}_{i=0}^d$ denote one of these decompositions.
For $0 \leq i \leq d$, the actions of $A$ and $A^*$ on $U_i$ are given in the table below; see \cite[Section~5]{IT2007RJ}.

\begin{equation*}\label{eq:action A,B on Ui}
\renewcommand{\arraystretch}{1.2} 
\begin{tabular}{c|c|c}
	decomposition name & action of $A$ on $U_i$ & action of $A^*$ on $U_i$ \\
	\hline \hline
	$[0, d]$ 		& $AU_i \subseteq U_i$ & $A^*U_i \subseteq U_{i-1} + U_i + U_{i+1}$ \\
	$[0^*, d^*]$ 	& $AU_i \subseteq U_{i-1} + U_i + U_{i+1}$ &  $A^*U_i \subseteq U_i$\\
	$[0^*, 0]$ 		& $AU_i \subseteq U_{i} + U_{i+1}$ & $A^*U_i \subseteq U_{i-1} + U_i$\\
	$[0^*, d]$ 		& $AU_i \subseteq U_{i} + U_{i+1}$ & $A^*U_i \subseteq U_{i-1} + U_i$ \\
	$[d^*, 0]$ 		& $AU_i \subseteq U_{i} + U_{i+1}$ & $A^*U_i \subseteq U_{i-1} + U_i$ \\
	$[d^*, d]$ 		& $AU_i \subseteq U_{i} + U_{i+1}$ & $A^*U_i \subseteq U_{i-1} + U_i$ 
\end{tabular}
\end{equation*}

\noindent
\begin{minipage}{\textwidth}

The following diagrams show how the actions described in the table above are represented on the tetrahedron diagram.

\begin{center}
\begin{tabular}{ccc}
\scalemath{0.7}{
\begin{tikzpicture}
  [scale=1,thick,auto=left, every node/.style={circle, inner sep=2pt}] 
  \node[fill=black,label=below:{\Large$E_0V$}] (b0) at (0,0) {};
  \node[fill=black,label=below:{}] (b1) at (1,0)  {};
  \node[fill=black,label=below:{}] (b2) at (2,0)  {};
  \node[fill=black,label=below:{}] (b3) at (3,0) {};
  \node[fill=black,label=below:{}] (b4) at (4,0)  {};
  \node[fill=black,label=below:{}] (b5) at (5,0)  {};
  \node[fill=black,label=below:{}] (b6) at (6,0)  {};
  \node[fill=black,label=below:{\Large$E_dV$}] (b7) at (7,0)  {};
  
  \node[fill=black,label=above:{\Large$E_0^*V$}] (t0) at (0,7)  {};
  \node[fill=black,label=above:{}] (t1) at (1,7)  {};
  \node[fill=black,label=above:{}] (t2) at (2,7)  {};
  \node[fill=black,label=above:{}] (t3) at (3,7)  {};
  \node[fill=black,label=above:{}] (t4) at (4,7)  {};
  \node[fill=black,label=above:{}] (t5) at (5,7)  {}; 
  \node[fill=black,label=above:{}] (t6) at (6,7)  {}; 
  \node[fill=black,label=above:{\Large$E_d^*V$}] (t7) at (7,7)  {}; 

  \node[fill=black,label=left:{}] (l1) at (0,1)  {};
  \node[fill=black,label=left:{}] (l2) at (0,2)  {};
  \node[fill=black,label=left:{}] (l3) at (0,3)  {};
  \node[fill=black,label=left:{}] (l4) at (0,4)  {};
  \node[fill=black,label=left:{}] (l5) at (0,5)  {};
  \node[fill=black,label=left:{}] (l6) at (0,6)  {};

  \node[fill=black,label=right:{}] (r1) at (7,1)  {};
  \node[fill=black,label=right:{}] (r2) at (7,2)  {};
  \node[fill=black,label=right:{}] (r3) at (7,3)  {};
  \node[fill=black,label=right:{}] (r4) at (7,4)  {};
  \node[fill=black,label=right:{}] (r5) at (7,5)  {};
  \node[fill=black,label=right:{}] (r6) at (7,6)  {};

  \node[fill=black,label=right:{}] (d1) at (1,1)  {};
  \node[fill=black,label=right:{}] (d2) at (1.9,1.9)  {};
  \node[fill=black,label=right:{}] (d3) at (2.8,2.8)  {};
  \node[fill=black,label=right:{}] (d4) at (4.2,4.2)  {};
  \node[fill=black,label=right:{}] (d5) at (5.1,5.1)  {};
  \node[fill=black,label=right:{}] (d6) at (6,6)  {};

  \node[fill=black,label=right:{}] (e1) at (1,6)  {};
  \node[fill=black,label=right:{}] (e2) at (1.9,5.1)  {};
  \node[fill=black,label=right:{}] (e3) at (2.8,4.2)  {};
  \node[fill=black,label=right:{}] (e4) at (4.2,2.8)  {};
  \node[fill=black,label=right:{}] (e5) at (5.1,1.9)  {};
  \node[fill=black,label=right:{}] (e6) at (6,1)  {};

  \foreach \from/\to in
  {b0/b1, b1/b2, b2/b3, b3/b4, b4/b5, b5/b6, b6/b7,
  t0/t1, t1/t2, t2/t3, t3/t4, t4/t5, t5/t6, t6/t7,                  
  b0/l1, l1/l2, l2/l3, l3/l4, l4/l5, l5/l6, l6/t0,
  t7/r1, r1/r2, r2/r3, r3/r4, r4/r5, r5/r6, r6/b7,
  t0/e1, e1/e2, e2/e3, e3/e4, e4/e5, e5/e6, e6/b7,
  b0/d1, d1/d2, d2/d3, d3/d4, d4/d5, d5/d6, d6/t7}       
      \draw (\from) -- (\to);        

     \draw (3,7) ellipse (0.2cm and 0.2cm);
     \draw (3,7) ellipse (0.3cm and 0.3cm);        
     \draw (2,7) ellipse (0.3cm and 0.3cm);
     \draw (4,7) ellipse (0.3cm and 0.3cm);
     
     \draw (3,0) ellipse (0.2cm and 0.2cm);
     \draw (3,0) ellipse (0.3cm and 0.3cm);
     
     \draw (0,5) ellipse (0.2cm and 0.2cm);
     \draw (0,5) ellipse (0.3cm and 0.3cm);        
     \draw (0,4) ellipse (0.3cm and 0.3cm);
     
     \draw (7,5) ellipse (0.2cm and 0.2cm);
     \draw (7,5) ellipse (0.3cm and 0.3cm);        
     \draw (7,4) ellipse (0.3cm and 0.3cm);
     
     \draw (1.9,5.1) ellipse (0.2cm and 0.2cm);
     \draw (1.9,5.1) ellipse (0.3cm and 0.3cm);        
     \draw (2.8,4.2) ellipse (0.3cm and 0.3cm);
     
     \draw (5.1,5.1) ellipse (0.2cm and 0.2cm);
     \draw (5.1,5.1) ellipse (0.3cm and 0.3cm);        
     \draw (4.2,4.2) ellipse (0.3cm and 0.3cm);
\end{tikzpicture}}

& \qquad \qquad &

\scalemath{0.7}{
\begin{tikzpicture}
  [scale=1,thick,auto=left, every node/.style={circle, inner sep=2pt}] 
  \node[fill=black,label=below:{\Large$E_0V$}] (b0) at (0,0) {};
  \node[fill=black,label=below:{}] (b1) at (1,0)  {};
  \node[fill=black,label=below:{}] (b2) at (2,0)  {};
  \node[fill=black,label=below:{}] (b3) at (3,0) {};
  \node[fill=black,label=below:{}] (b4) at (4,0)  {};
  \node[fill=black,label=below:{}] (b5) at (5,0)  {};
  \node[fill=black,label=below:{}] (b6) at (6,0)  {};
  \node[fill=black,label=below:{\Large$E_dV$}] (b7) at (7,0)  {};
  
  \node[fill=black,label=above:{\Large$E_0^*V$}] (t0) at (0,7)  {};
  \node[fill=black,label=above:{}] (t1) at (1,7)  {};
  \node[fill=black,label=above:{}] (t2) at (2,7)  {};
  \node[fill=black,label=above:{}] (t3) at (3,7)  {};
  \node[fill=black,label=above:{}] (t4) at (4,7)  {};
  \node[fill=black,label=above:{}] (t5) at (5,7)  {}; 
  \node[fill=black,label=above:{}] (t6) at (6,7)  {}; 
  \node[fill=black,label=above:{\Large$E_d^*V$}] (t7) at (7,7)  {}; 

  \node[fill=black,label=left:{}] (l1) at (0,1)  {};
  \node[fill=black,label=left:{}] (l2) at (0,2)  {};
  \node[fill=black,label=left:{}] (l3) at (0,3)  {};
  \node[fill=black,label=left:{}] (l4) at (0,4)  {};
  \node[fill=black,label=left:{}] (l5) at (0,5)  {};
  \node[fill=black,label=left:{}] (l6) at (0,6)  {};

  \node[fill=black,label=right:{}] (r1) at (7,1)  {};
  \node[fill=black,label=right:{}] (r2) at (7,2)  {};
  \node[fill=black,label=right:{}] (r3) at (7,3)  {};
  \node[fill=black,label=right:{}] (r4) at (7,4)  {};
  \node[fill=black,label=right:{}] (r5) at (7,5)  {};
  \node[fill=black,label=right:{}] (r6) at (7,6)  {};

  \node[fill=black,label=right:{}] (d1) at (1,1)  {};
  \node[fill=black,label=right:{}] (d2) at (1.9,1.9)  {};
  \node[fill=black,label=right:{}] (d3) at (2.8,2.8)  {};
  \node[fill=black,label=right:{}] (d4) at (4.2,4.2)  {};
  \node[fill=black,label=right:{}] (d5) at (5.1,5.1)  {};
  \node[fill=black,label=right:{}] (d6) at (6,6)  {};

  \node[fill=black,label=right:{}] (e1) at (1,6)  {};
  \node[fill=black,label=right:{}] (e2) at (1.9,5.1)  {};
  \node[fill=black,label=right:{}] (e3) at (2.8,4.2)  {};
  \node[fill=black,label=right:{}] (e4) at (4.2,2.8)  {};
  \node[fill=black,label=right:{}] (e5) at (5.1,1.9)  {};
  \node[fill=black,label=right:{}] (e6) at (6,1)  {};

  \foreach \from/\to in
  {b0/b1, b1/b2, b2/b3, b3/b4, b4/b5, b5/b6, b6/b7,
  t0/t1, t1/t2, t2/t3, t3/t4, t4/t5, t5/t6, t6/t7,                  
  b0/l1, l1/l2, l2/l3, l3/l4, l4/l5, l5/l6, l6/t0,
  t7/r1, r1/r2, r2/r3, r3/r4, r4/r5, r5/r6, r6/b7,
  t0/e1, e1/e2, e2/e3, e3/e4, e4/e5, e5/e6, e6/b7,
  b0/d1, d1/d2, d2/d3, d3/d4, d4/d5, d5/d6, d6/t7}       
      \draw (\from) -- (\to);        

     \draw (4,7) ellipse (0.2cm and 0.2cm);
     \draw (4,7) ellipse (0.3cm and 0.3cm);
          
     \draw (4,0) ellipse (0.2cm and 0.2cm);
     \draw (4,0) ellipse (0.3cm and 0.3cm);        
     \draw (3,0) ellipse (0.3cm and 0.3cm);
     \draw (5,0) ellipse (0.3cm and 0.3cm);
     
     \draw (0,2) ellipse (0.2cm and 0.2cm);
     \draw (0,2) ellipse (0.3cm and 0.3cm);        
     \draw (0,3) ellipse (0.3cm and 0.3cm);
     
     \draw (7,2) ellipse (0.2cm and 0.2cm);
     \draw (7,2) ellipse (0.3cm and 0.3cm);        
     \draw (7,3) ellipse (0.3cm and 0.3cm);
     
     \draw (1.9,1.9) ellipse (0.2cm and 0.2cm);
     \draw (1.9,1.9) ellipse (0.3cm and 0.3cm);        
     \draw (2.8,2.8) ellipse (0.3cm and 0.3cm);
     
     \draw (5.1,1.9) ellipse (0.2cm and 0.2cm);
     \draw (5.1,1.9) ellipse (0.3cm and 0.3cm);
     \draw (4.2,2.8) ellipse (0.3cm and 0.3cm);        

\end{tikzpicture}} \\
Action of $A$ & \qquad \qquad & Action of $A^*$
\end{tabular}
\end{center}
In the left diagram, a small circle marks a component $U_i$ on which $A$ acts.
The large circles mark the components whose sum contains $AU_i$.
The right diagram similarly describes the action of $A^*$.

\end{minipage}

\section{TD systems and $\mathcal{O}$-modules}\label{sec:TDsystemO-mod}

In this section, we return our attention to the Onsager Lie algebra $\mathcal{O}$.
We use the facts about TD systems reviewed in Sections~\ref{Sec:TD_pair}--\ref{Sec:decomp-diag} to describe the finite-dimensional irreducible $\mathcal{O}$-modules.
For the rest of the paper, assume that $\mathbb{F}$ is algebraically closed with characteristic zero.

\begin{definition}\label{def:arithmetic2}
Let $\theta_0, \theta_1, \ldots, \theta_d$ denote a finite sequence of scalars in $\mathbb{F}$.
We say that this sequence is \emph{arithmetic with common difference $2$} whenever $\theta_{i-1} - \theta_i = 2$ for $1 \leq i \leq d$.
\end{definition}

\begin{lemma}[{\cite[Theorem~2.4]{2007Har}}]\label{lem:Thm2.4Hartwig}
Let $V$ denote a nonzero finite-dimensional vector space, and let $A,A^* \in \End(V)$.
Then the following {\rm(i)}, {\rm(ii)} are equivalent.
\begin{itemize}
	\item[\rm(i)] $A,A^*$ is a TD pair on $V$ with an eigenvalue sequence and a dual eigenvalue sequence that are both arithmetic with common difference $2$.
	\item[\rm(ii)] There exists an irreducible $\mathcal{O}$-module structure on $V$ such that the first and second standard $\mathcal{O}$-generators act on $V$ as $A$ and $A^*$, respectively.
\end{itemize}
\end{lemma}

Let $A,A^*$ denote a TD pair on $V$ of diameter $d$.
Assume that the TD pair $A, A^*$ has an eigenvalue sequence and a dual eigenvalue sequence that are both arithmetic with common difference $2$.
By the discussion in the last paragraph of Section~\ref{Sec:TD_pair}, we may choose a TD system
$$
    \Phi=(A;\{E_i\}_{i=0}^d;A^*;\{E_i^*\}_{i=0}^d)
$$
with eigenvalue sequence $\{\theta_i\}_{i=0}^d$ and dual eigenvalue sequence $\{\theta_i^*\}_{i=0}^d$ of the form
\begin{equation*}
    \theta_i=2i-d+\alpha,
    \qquad
    \theta_i^*=d-2i+\alpha^*
    \qquad (0\leq i\leq d),
\end{equation*}
where $\alpha,\alpha^*\in\mathbb{F}$.
Observe that
\begin{equation*}
    \theta_i-\theta_j=2(i-j),
    \qquad
    \theta_i^*-\theta_j^*=2(j-i)
    \qquad (0\leq i,j\leq d).
\end{equation*}
By Lemma~\ref{lem:Thm2.4Hartwig}, we turn $V$ into an irreducible $\mathcal{O}$-module on which the first and second standard $\mathcal{O}$-generators $A,B$ act as $A,A^*$, respectively.

\begin{lemma}\label{lem:triple rels}
For $k \in \mathbb{N}$ and $0 \leq i,j \leq d$, the following {\rm(i)}--{\rm(iv)} hold for the $\mathcal{O}$-module $V$:
\begin{itemize}
\item[\rm(i)] $A_k$ satisfies
	\begin{align}
		(i - j) E_iA_kE_j & = 0, \label{eq(1):t.rels A_k(i,j)}\\
		(i - j)(i - j - 1)(i - j + 1) E_i^*A_kE_j^* & = 0. \label{eq(2):t.rels A_k(i,j)}
	\end{align}
\item[\rm(ii)] $B_k$ satisfies
	\begin{align}
		(i - j)(i - j - 1)(i - j + 1) E_iB_kE_j & = 0, \label{eq(1):t.rels B_k(i,j)}\\
		(i - j) E^*_iB_kE_j^* & = 0. \label{eq(2):t.rels B_k(i,j)}
	\end{align}
\item[\rm(iii)] $\psi_{k+1}$ satisfies
	\begin{align}
		(i - j)(i - j + 1) E_i \psi_{k+1} E_j & = 0, \label{eq(1):t.rels psi_k(i,j)}\\
		(i - j)(i - j + 1) E_i^* \psi_{k+1} E_j^* & = 0. \label{eq(2):t.rels psi_k(i,j)}
	\end{align}	
\item[\rm(iv)] $\tilde\psi_{k+1}$ satisfies
	\begin{align}
		(i - j)(i - j - 1) E_i \tilde\psi_{k+1} E_j & = 0, \label{eq(1):t.rels tilde psi_k(i,j)}\\
		(i - j)(i - j - 1) E_i^* \tilde\psi_{k+1} E_j^* & = 0. \label{eq(2):t.rels tilde psi_k(i,j)}
	\end{align}	
\end{itemize}
\end{lemma}
\begin{proof}
For each equation in Lemma \ref{lem:some rels A,B} that involves $A$, multiply both sides  on the left by $E_i$ and on the right by $E_j$.
Simplify the result using $E_iA = \theta_i E_i$ and $AE_j = \theta_j E_j$.
Evaluate the result further using $\theta_i - \theta_j = 2(i - j)$ to get the equations \eqref{eq(1):t.rels A_k(i,j)}, \eqref{eq(1):t.rels B_k(i,j)}, \eqref{eq(1):t.rels psi_k(i,j)}, \eqref{eq(1):t.rels tilde psi_k(i,j)}.
Similarly, for each equation in Lemma \ref{lem:some rels A,B} that involves $B$, multiply both sides on the left by $E^*_i$ and on the right by $E_j^*$.
Simplify the result using $E_i^*B = \theta_i^* E_i^*$ and $BE_j^* = \theta_j^* E_j^*$.
Evaluate the result further using $\theta_i^* - \theta_j^* = 2(j - i)$ to get the equations \eqref{eq(2):t.rels A_k(i,j)}, \eqref{eq(2):t.rels B_k(i,j)}, \eqref{eq(2):t.rels psi_k(i,j)}, \eqref{eq(2):t.rels tilde psi_k(i,j)}.
The result follows.
\end{proof}

\begin{lemma}\label{lem:triple rels(2)}
For $k \in \mathbb{N}$ and $0 \leq i,j \leq d$, the following {\rm(i)}--{\rm(iv)} hold for the $\mathcal{O}$-module $V$:
\begin{itemize}
\item[\rm(i)] $A_k$ satisfies
	\begin{align}
		E_iA_kE_j & = 0 \qquad \text{if} \quad i \ne j, \label{eq(1):triple rels(2)} \\
		E_i^*A_kE_j^* & = 0 \qquad \text{if} \quad |i - j| > 1. \label{eq(2):triple rels(2)}
	\end{align}
\item[\rm(ii)] $B_k$ satisfies
	\begin{align*}
		E_iB_kE_j & = 0 \qquad \text{if} \quad |i - j| > 1, \\
		E^*_iB_kE_j^* & = 0 \qquad \text{if} \quad i \ne j.
	\end{align*}
\item[\rm(iii)] $\psi_{k+1}$ satisfies
	\begin{align*}
		E_i \psi_{k+1} E_j & = 0 \qquad \text{if} \quad j - i \notin \{0, 1\}, \\
		E_i^* \psi_{k+1} E_j^* & = 0 \qquad \text{if} \quad j - i \notin \{0, 1\}.
	\end{align*}	
\item[\rm(iv)] $\tilde\psi_{k+1}$ satisfies
	\begin{align*}
		E_i \tilde\psi_{k+1} E_j & = 0 \qquad \text{if} \quad i - j \notin \{0, 1\}, \\
		E_i^* \tilde\psi_{k+1} E_j^* & = 0 \qquad \text{if} \quad i - j \notin \{0, 1\}.
	\end{align*}	
\end{itemize}
\end{lemma}
\begin{proof}
By Lemma~\ref{lem:triple rels}.
\end{proof}

For notational convenience, set $E_{-1}=E_{d+1}=E_{-1}^*=E_{d+1}^*=0$.

\begin{lemma}\label{lem:act.alt.elts}
For $k \in \mathbb{N}$ and $0 \leq j \leq d$, the following {\rm(i)}--{\rm(iv)} hold for the $\mathcal{O}$-module $V$:

\begin{itemize}
\item[\rm(i)] $A_k$ satisfies
	\begin{align*}
		A_kE_jV & \subseteq E_jV, \\
		A_kE_j^*V & \subseteq E_{j-1}^*V + E_j^*V + E_{j+1}^*V.
	\end{align*}
\item[\rm(ii)] $B_k$ satisfies
	\begin{align*}
		B_kE_jV & \subseteq E_{j-1}V + E_jV + E_{j+1}V, \\
		B_kE_j^*V & \subseteq E_j^*V.
	\end{align*}
\item[\rm(iii)] $\psi_{k+1}$ satisfies
	\begin{align*}
		\psi_{k+1} E_jV & \subseteq E_{j-1}V + E_jV, \\
		\psi_{k+1} E_j^*V & \subseteq E_{j-1}^*V + E_j^*V.
	\end{align*}	
\item[\rm(iv)] $\tilde\psi_{k+1}$ satisfies
	\begin{align*}
		\tilde\psi_{k+1} E_j V & \subseteq E_j V + E_{j+1} V, \\
		\tilde\psi_{k+1} E_j^* V& \subseteq E_j^* V + E_{j+1}^* V.
	\end{align*}	
\end{itemize}
\end{lemma}
\begin{proof}
(i) By \eqref{eq(1):triple rels(2)}, we obtain
$$
	A_k E_j V = I A_k E_j V \subseteq \sum_{i=0}^d E_i A_k E_j V = E_j A_k E_j V \subseteq E_jV.
$$
Next, by \eqref{eq(2):triple rels(2)}
$$
	A_k E_j^* V = I A_k E_j^* V \subseteq \sum_{i=0}^d E_i^* A_k E_j^* V = \sum_{i=j-1}^{j+1} E_i^* A_k E_j^* V \subseteq E_{j-1}^*V + E_{j}^*V + E_{j+1}^*V.
$$
(ii)--(iv) Similar to (i) using Lemma~\ref{lem:triple rels(2)}(ii)--(iv).
\end{proof}

We finish this section by describing how the elements in \eqref{eq:module elements} act on the flags $[0]$, $[d]$, $[0^*]$, $[d^*]$ described in \eqref{table:4flags}.

\begin{lemma}\label{lem:flag behavior}
For $k\in\mathbb{N}$, the actions of $A_k$, $B_k$, $\psi_{k+1}$, $\tilde\psi_{k+1}$ on the flags $[0]$, $[d]$, $[0^*]$, $[d^*]$ are described in the table below.
\begin{equation*}
\renewcommand{\arraystretch}{1.2}
\begin{tabular}{c|c|c|c|c}
	element & $[0]$ & $[d]$ & $[0^*]$ & $[d^*]$ \\
	\hline\hline
	$A_k$ & stabilizes & stabilizes & raises & raises \\
	$B_k$ & raises & raises & stabilizes & stabilizes \\
	$\psi_{k+1}$ & stabilizes & raises & stabilizes & raises \\
	$\tilde\psi_{k+1}$ & raises & stabilizes & raises & stabilizes
\end{tabular}
\end{equation*}

\end{lemma}

\begin{proof}
By \eqref{table:4flags}, the pairs of opposite flags $[0],[d]$ and $[0^*],[d^*]$ arise from the decompositions $\{E_iV\}_{i=0}^d$ and $\{E_i^*V\}_{i=0}^d$, respectively.
Apply Lemma~\ref{lem:raise,stabilizes} to these pairs using Lemma~\ref{lem:act.alt.elts}.
\end{proof}

\section{The tetrahedron diagram revisited}\label{Sec:action-on-V}

We continue to discuss the finite-dimensional irreducible $\mathcal{O}$-module $V$ from Section \ref{sec:TDsystemO-mod}.
In this section, we describe how the elements in \eqref{eq:module elements} act on the six decompositions of $V$ described in \eqref{table:6decomp}.
We use the tetrahedron diagram to illustrate these actions.
Let $\{U_i\}_{i=0}^d$ be one of the decompositions of $V$ listed in \eqref{table:6decomp}, and let $[u,v]$ denote the name of this decomposition $\{U_i\}_{i=0}^d$.
By construction, the flag $[u]$ is induced by $\{U_i\}_{i=0}^d$ and the flag $[v]$ is induced by $\{U_{d-i}\}_{i=0}^d$.

\smallskip
For $k \in \mathbb{N}$, we now describe the actions of $A_k$ and $B_k$ on $\{U_i\}_{i=0}^d$.

\begin{theorem}\label{thm:action A_k B_k on Ui}
For $k\in\mathbb{N}$ and $0\leq i\leq d$, the actions of $A_k$ and $B_k$ on $U_i$ are described in the table below.

\begin{equation*}
\renewcommand{\arraystretch}{1.2}
\begin{tabular}{c|c|c}

	decomposition name
	& action of $A_k$ on $U_i$
	& action of $B_k$ on $U_i$ \\

	\hline \hline

	$[0,d]$
	& $A_kU_i\subseteq U_i$
	& $B_kU_i\subseteq U_{i-1}+U_i+U_{i+1}$ \\

	$[0^*,d^*]$
	& $A_kU_i\subseteq U_{i-1}+U_i+U_{i+1}$
	& $B_kU_i\subseteq U_i$ \\

	$[0^*,0]$
	& $A_kU_i\subseteq U_i+U_{i+1}$
	& $B_kU_i\subseteq U_{i-1}+U_i$ \\

	$[0^*,d]$
	& $A_kU_i\subseteq U_i+U_{i+1}$
	& $B_kU_i\subseteq U_{i-1}+U_i$ \\

	$[d^*,0]$
	& $A_kU_i\subseteq U_i+U_{i+1}$
	& $B_kU_i\subseteq U_{i-1}+U_i$ \\

	$[d^*,d]$
	& $A_kU_i\subseteq U_i+U_{i+1}$
	& $B_kU_i\subseteq U_{i-1}+U_i$

\end{tabular}
\end{equation*}
\end{theorem}

\begin{proof}
Apply Lemma~\ref{lem:raise,stabilizes} using Lemma~\ref{lem:flag behavior}.
\end{proof}

\noindent
\begin{minipage}{\textwidth}

In the figure below, we use the tetrahedron diagram to illustrate Theorem~\ref{thm:action A_k B_k on Ui}.

\begin{center}
\begin{tabular}{ccc}
\scalemath{0.7}{
\begin{tikzpicture}
  [scale=1,thick,auto=left, every node/.style={circle, inner sep=2pt}] 
  \node[fill=black,label=below:{\Large$E_0V$}] (b0) at (0,0) {};
  \node[fill=black,label=below:{}] (b1) at (1,0)  {};
  \node[fill=black,label=below:{}] (b2) at (2,0)  {};
  \node[fill=black,label=below:{}] (b3) at (3,0) {};
  \node[fill=black,label=below:{}] (b4) at (4,0)  {};
  \node[fill=black,label=below:{}] (b5) at (5,0)  {};
  \node[fill=black,label=below:{}] (b6) at (6,0)  {};
  \node[fill=black,label=below:{\Large$E_dV$}] (b7) at (7,0)  {};
  
  \node[fill=black,label=above:{\Large$E_0^*V$}] (t0) at (0,7)  {};
  \node[fill=black,label=above:{}] (t1) at (1,7)  {};
  \node[fill=black,label=above:{}] (t2) at (2,7)  {};
  \node[fill=black,label=above:{}] (t3) at (3,7)  {};
  \node[fill=black,label=above:{}] (t4) at (4,7)  {};
  \node[fill=black,label=above:{}] (t5) at (5,7)  {}; 
  \node[fill=black,label=above:{}] (t6) at (6,7)  {}; 
  \node[fill=black,label=above:{\Large$E_d^*V$}] (t7) at (7,7)  {}; 

  \node[fill=black,label=left:{}] (l1) at (0,1)  {};
  \node[fill=black,label=left:{}] (l2) at (0,2)  {};
  \node[fill=black,label=left:{}] (l3) at (0,3)  {};
  \node[fill=black,label=left:{}] (l4) at (0,4)  {};
  \node[fill=black,label=left:{}] (l5) at (0,5)  {};
  \node[fill=black,label=left:{}] (l6) at (0,6)  {};

  \node[fill=black,label=right:{}] (r1) at (7,1)  {};
  \node[fill=black,label=right:{}] (r2) at (7,2)  {};
  \node[fill=black,label=right:{}] (r3) at (7,3)  {};
  \node[fill=black,label=right:{}] (r4) at (7,4)  {};
  \node[fill=black,label=right:{}] (r5) at (7,5)  {};
  \node[fill=black,label=right:{}] (r6) at (7,6)  {};

  \node[fill=black,label=right:{}] (d1) at (1,1)  {};
  \node[fill=black,label=right:{}] (d2) at (1.9,1.9)  {};
  \node[fill=black,label=right:{}] (d3) at (2.8,2.8)  {};
  \node[fill=black,label=right:{}] (d4) at (4.2,4.2)  {};
  \node[fill=black,label=right:{}] (d5) at (5.1,5.1)  {};
  \node[fill=black,label=right:{}] (d6) at (6,6)  {};

  \node[fill=black,label=right:{}] (e1) at (1,6)  {};
  \node[fill=black,label=right:{}] (e2) at (1.9,5.1)  {};
  \node[fill=black,label=right:{}] (e3) at (2.8,4.2)  {};
  \node[fill=black,label=right:{}] (e4) at (4.2,2.8)  {};
  \node[fill=black,label=right:{}] (e5) at (5.1,1.9)  {};
  \node[fill=black,label=right:{}] (e6) at (6,1)  {};

  \foreach \from/\to in
  {b0/b1, b1/b2, b2/b3, b3/b4, b4/b5, b5/b6, b6/b7,
  t0/t1, t1/t2, t2/t3, t3/t4, t4/t5, t5/t6, t6/t7,                  
  b0/l1, l1/l2, l2/l3, l3/l4, l4/l5, l5/l6, l6/t0,
  t7/r1, r1/r2, r2/r3, r3/r4, r4/r5, r5/r6, r6/b7,
  t0/e1, e1/e2, e2/e3, e3/e4, e4/e5, e5/e6, e6/b7,
  b0/d1, d1/d2, d2/d3, d3/d4, d4/d5, d5/d6, d6/t7}       
      \draw (\from) -- (\to);

     \draw (3,7) ellipse (0.2cm and 0.2cm);
     \draw (3,7) ellipse (0.3cm and 0.3cm);        
     \draw (2,7) ellipse (0.3cm and 0.3cm);
     \draw (4,7) ellipse (0.3cm and 0.3cm);
     
     \draw (3,0) ellipse (0.2cm and 0.2cm);
     \draw (3,0) ellipse (0.3cm and 0.3cm);
     
     \draw (0,5) ellipse (0.2cm and 0.2cm);
     \draw (0,5) ellipse (0.3cm and 0.3cm);        
     \draw (0,4) ellipse (0.3cm and 0.3cm);
     
     \draw (7,5) ellipse (0.2cm and 0.2cm);
     \draw (7,5) ellipse (0.3cm and 0.3cm);        
     \draw (7,4) ellipse (0.3cm and 0.3cm);
     
     \draw (1.9,5.1) ellipse (0.2cm and 0.2cm);
     \draw (1.9,5.1) ellipse (0.3cm and 0.3cm);        
     \draw (2.8,4.2) ellipse (0.3cm and 0.3cm);
     
     \draw (5.1,5.1) ellipse (0.2cm and 0.2cm);
     \draw (5.1,5.1) ellipse (0.3cm and 0.3cm);        
     \draw (4.2,4.2) ellipse (0.3cm and 0.3cm);
\end{tikzpicture}}

& \qquad \qquad &

\scalemath{0.7}{
\begin{tikzpicture}
  [scale=1,thick,auto=left, every node/.style={circle, inner sep=2pt}] 
  \node[fill=black,label=below:{\Large$E_0V$}] (b0) at (0,0) {};
  \node[fill=black,label=below:{}] (b1) at (1,0)  {};
  \node[fill=black,label=below:{}] (b2) at (2,0)  {};
  \node[fill=black,label=below:{}] (b3) at (3,0) {};
  \node[fill=black,label=below:{}] (b4) at (4,0)  {};
  \node[fill=black,label=below:{}] (b5) at (5,0)  {};
  \node[fill=black,label=below:{}] (b6) at (6,0)  {};
  \node[fill=black,label=below:{\Large$E_dV$}] (b7) at (7,0)  {};
  
  \node[fill=black,label=above:{\Large$E_0^*V$}] (t0) at (0,7)  {};
  \node[fill=black,label=above:{}] (t1) at (1,7)  {};
  \node[fill=black,label=above:{}] (t2) at (2,7)  {};
  \node[fill=black,label=above:{}] (t3) at (3,7)  {};
  \node[fill=black,label=above:{}] (t4) at (4,7)  {};
  \node[fill=black,label=above:{}] (t5) at (5,7)  {}; 
  \node[fill=black,label=above:{}] (t6) at (6,7)  {}; 
  \node[fill=black,label=above:{\Large$E_d^*V$}] (t7) at (7,7)  {}; 

  \node[fill=black,label=left:{}] (l1) at (0,1)  {};
  \node[fill=black,label=left:{}] (l2) at (0,2)  {};
  \node[fill=black,label=left:{}] (l3) at (0,3)  {};
  \node[fill=black,label=left:{}] (l4) at (0,4)  {};
  \node[fill=black,label=left:{}] (l5) at (0,5)  {};
  \node[fill=black,label=left:{}] (l6) at (0,6)  {};

  \node[fill=black,label=right:{}] (r1) at (7,1)  {};
  \node[fill=black,label=right:{}] (r2) at (7,2)  {};
  \node[fill=black,label=right:{}] (r3) at (7,3)  {};
  \node[fill=black,label=right:{}] (r4) at (7,4)  {};
  \node[fill=black,label=right:{}] (r5) at (7,5)  {};
  \node[fill=black,label=right:{}] (r6) at (7,6)  {};

  \node[fill=black,label=right:{}] (d1) at (1,1)  {};
  \node[fill=black,label=right:{}] (d2) at (1.9,1.9)  {};
  \node[fill=black,label=right:{}] (d3) at (2.8,2.8)  {};
  \node[fill=black,label=right:{}] (d4) at (4.2,4.2)  {};
  \node[fill=black,label=right:{}] (d5) at (5.1,5.1)  {};
  \node[fill=black,label=right:{}] (d6) at (6,6)  {};

  \node[fill=black,label=right:{}] (e1) at (1,6)  {};
  \node[fill=black,label=right:{}] (e2) at (1.9,5.1)  {};
  \node[fill=black,label=right:{}] (e3) at (2.8,4.2)  {};
  \node[fill=black,label=right:{}] (e4) at (4.2,2.8)  {};
  \node[fill=black,label=right:{}] (e5) at (5.1,1.9)  {};
  \node[fill=black,label=right:{}] (e6) at (6,1)  {};

  \foreach \from/\to in
  {b0/b1, b1/b2, b2/b3, b3/b4, b4/b5, b5/b6, b6/b7,
  t0/t1, t1/t2, t2/t3, t3/t4, t4/t5, t5/t6, t6/t7,                  
  b0/l1, l1/l2, l2/l3, l3/l4, l4/l5, l5/l6, l6/t0,
  t7/r1, r1/r2, r2/r3, r3/r4, r4/r5, r5/r6, r6/b7,
  t0/e1, e1/e2, e2/e3, e3/e4, e4/e5, e5/e6, e6/b7,
  b0/d1, d1/d2, d2/d3, d3/d4, d4/d5, d5/d6, d6/t7}       
      \draw (\from) -- (\to);

     \draw (4,7) ellipse (0.2cm and 0.2cm);
     \draw (4,7) ellipse (0.3cm and 0.3cm);
          
     \draw (4,0) ellipse (0.2cm and 0.2cm);
     \draw (4,0) ellipse (0.3cm and 0.3cm);        
     \draw (3,0) ellipse (0.3cm and 0.3cm);
     \draw (5,0) ellipse (0.3cm and 0.3cm);
     
     \draw (0,2) ellipse (0.2cm and 0.2cm);
     \draw (0,2) ellipse (0.3cm and 0.3cm);        
     \draw (0,3) ellipse (0.3cm and 0.3cm);
     
     \draw (7,2) ellipse (0.2cm and 0.2cm);
     \draw (7,2) ellipse (0.3cm and 0.3cm);        
     \draw (7,3) ellipse (0.3cm and 0.3cm);
     
     \draw (1.9,1.9) ellipse (0.2cm and 0.2cm);
     \draw (1.9,1.9) ellipse (0.3cm and 0.3cm);        
     \draw (2.8,2.8) ellipse (0.3cm and 0.3cm);
     
     \draw (5.1,1.9) ellipse (0.2cm and 0.2cm);
     \draw (5.1,1.9) ellipse (0.3cm and 0.3cm);
     \draw (4.2,2.8) ellipse (0.3cm and 0.3cm);        

\end{tikzpicture}} \\
Action of $A_k$ & \qquad \qquad & Action of $B_k$
\end{tabular}
\end{center}

\noindent
In the left diagram, a small circle marks a component $U_i$ on which $A_k$ acts.
The large circles mark the components whose sum contains $A_kU_i$.
The right diagram similarly describes the action of $B_k$.

\end{minipage}

\smallskip
For $k \in \mathbb{N}$ we now describe the actions of $\psi_{k+1}$ and $\tilde\psi_{k+1}$.

\begin{theorem}\label{thm:action psi_k tilde.psi_k on Ui}
For $k\in\mathbb{N}$ and $0\leq i\leq d$, the actions of $\psi_{k+1}$ and $\tilde\psi_{k+1}$ on $U_i$ are described in the table below.

\begin{equation*}
\renewcommand{\arraystretch}{1.2}
\begin{tabular}{c|c|c}

	decomposition name
	& action of $\psi_{k+1}$ on $U_i$
	& action of $\tilde\psi_{k+1}$ on $U_i$ \\

	\hline \hline

	$[0,d]$
	& $\psi_{k+1}U_i\subseteq U_{i-1}+U_i$
	& $\tilde\psi_{k+1}U_i\subseteq U_i+U_{i+1}$ \\

	$[0^*,d^*]$
	& $\psi_{k+1}U_i\subseteq U_{i-1}+U_i$
	& $\tilde\psi_{k+1}U_i\subseteq U_i+U_{i+1}$ \\

	$[0^*,0]$
	& $\psi_{k+1}U_i\subseteq U_i$
	& $\tilde\psi_{k+1}U_i\subseteq U_{i-1}+U_i+U_{i+1}$ \\

	$[0^*,d]$
	& $\psi_{k+1}U_i\subseteq U_{i-1}+U_i$
	& $\tilde\psi_{k+1}U_i\subseteq U_i+U_{i+1}$ \\

	$[d^*,0]$
	& $\psi_{k+1}U_i\subseteq U_i+U_{i+1}$
	& $\tilde\psi_{k+1}U_i\subseteq U_{i-1}+U_i$ \\

	$[d^*,d]$
	& $\psi_{k+1}U_i\subseteq U_{i-1}+U_i+U_{i+1}$
	& $\tilde\psi_{k+1}U_i\subseteq U_i$

\end{tabular}
\end{equation*}
\end{theorem}

\begin{proof}
Apply Lemma~\ref{lem:raise,stabilizes} using Lemma~\ref{lem:flag behavior}.
\end{proof}

\noindent
\begin{minipage}{\textwidth}

In the figure below, we use the tetrahedron diagram to illustrate Theorem~\ref{thm:action psi_k tilde.psi_k on Ui}.

\begin{center}
\begin{tabular}{ccc}
\scalemath{0.7}{
\begin{tikzpicture}
  [scale=1,thick,auto=left, every node/.style={circle, inner sep=2pt}] 
  \node[fill=black,label=below:{\Large$E_0V$}] (b0) at (0,0) {};
  \node[fill=black,label=below:{}] (b1) at (1,0)  {};
  \node[fill=black,label=below:{}] (b2) at (2,0)  {};
  \node[fill=black,label=below:{}] (b3) at (3,0) {};
  \node[fill=black,label=below:{}] (b4) at (4,0)  {};
  \node[fill=black,label=below:{}] (b5) at (5,0)  {};
  \node[fill=black,label=below:{}] (b6) at (6,0)  {};
  \node[fill=black,label=below:{\Large$E_dV$}] (b7) at (7,0)  {};
  
  \node[fill=black,label=above:{\Large$E_0^*V$}] (t0) at (0,7)  {};
  \node[fill=black,label=above:{}] (t1) at (1,7)  {};
  \node[fill=black,label=above:{}] (t2) at (2,7)  {};
  \node[fill=black,label=above:{}] (t3) at (3,7)  {};
  \node[fill=black,label=above:{}] (t4) at (4,7)  {};
  \node[fill=black,label=above:{}] (t5) at (5,7)  {}; 
  \node[fill=black,label=above:{}] (t6) at (6,7)  {}; 
  \node[fill=black,label=above:{\Large$E_d^*V$}] (t7) at (7,7)  {}; 

  \node[fill=black,label=left:{}] (l1) at (0,1)  {};
  \node[fill=black,label=left:{}] (l2) at (0,2)  {};
  \node[fill=black,label=left:{}] (l3) at (0,3)  {};
  \node[fill=black,label=left:{}] (l4) at (0,4)  {};
  \node[fill=black,label=left:{}] (l5) at (0,5)  {};
  \node[fill=black,label=left:{}] (l6) at (0,6)  {};

  \node[fill=black,label=right:{}] (r1) at (7,1)  {};
  \node[fill=black,label=right:{}] (r2) at (7,2)  {};
  \node[fill=black,label=right:{}] (r3) at (7,3)  {};
  \node[fill=black,label=right:{}] (r4) at (7,4)  {};
  \node[fill=black,label=right:{}] (r5) at (7,5)  {};
  \node[fill=black,label=right:{}] (r6) at (7,6)  {};

  \node[fill=black,label=right:{}] (d1) at (1,1)  {};
  \node[fill=black,label=right:{}] (d2) at (1.9,1.9)  {};
  \node[fill=black,label=right:{}] (d3) at (2.8,2.8)  {};
  \node[fill=black,label=right:{}] (d4) at (4.2,4.2)  {};
  \node[fill=black,label=right:{}] (d5) at (5.1,5.1)  {};
  \node[fill=black,label=right:{}] (d6) at (6,6)  {};

  \node[fill=black,label=right:{}] (e1) at (1,6)  {};
  \node[fill=black,label=right:{}] (e2) at (1.9,5.1)  {};
  \node[fill=black,label=right:{}] (e3) at (2.8,4.2)  {};
  \node[fill=black,label=right:{}] (e4) at (4.2,2.8)  {};
  \node[fill=black,label=right:{}] (e5) at (5.1,1.9)  {};
  \node[fill=black,label=right:{}] (e6) at (6,1)  {};

  \foreach \from/\to in
  {b0/b1, b1/b2, b2/b3, b3/b4, b4/b5, b5/b6, b6/b7,
  t0/t1, t1/t2, t2/t3, t3/t4, t4/t5, t5/t6, t6/t7,                  
  b0/l1, l1/l2, l2/l3, l3/l4, l4/l5, l5/l6, l6/t0,
  t7/r1, r1/r2, r2/r3, r3/r4, r4/r5, r5/r6, r6/b7,
  t0/e1, e1/e2, e2/e3, e3/e4, e4/e5, e5/e6, e6/b7,
  b0/d1, d1/d2, d2/d3, d3/d4, d4/d5, d5/d6, d6/t7}       
      \draw (\from) -- (\to);

     \draw (3,7) ellipse (0.2cm and 0.2cm);
     \draw (3,7) ellipse (0.3cm and 0.3cm);        
     \draw (2,7) ellipse (0.3cm and 0.3cm);
     
     \draw (3,0) ellipse (0.2cm and 0.2cm);
     \draw (3,0) ellipse (0.3cm and 0.3cm);
     \draw (2,0) ellipse (0.3cm and 0.3cm);
     
     \draw (0,3) ellipse (0.2cm and 0.2cm);
     \draw (0,3) ellipse (0.3cm and 0.3cm);
     
     \draw (7,3) ellipse (0.2cm and 0.2cm);
     \draw (7,3) ellipse (0.3cm and 0.3cm);        
     \draw (7,4) ellipse (0.3cm and 0.3cm);
     \draw (7,2) ellipse (0.3cm and 0.3cm);
     
     \draw (5.1,1.9) ellipse (0.2cm and 0.2cm);
     \draw (4.2,2.8) ellipse (0.3cm and 0.3cm);        
     \draw (5.1,1.9) ellipse (0.3cm and 0.3cm);
     
     \draw (5.1,5.1) ellipse (0.2cm and 0.2cm);
     \draw (5.1,5.1) ellipse (0.3cm and 0.3cm);        
     \draw (4.2,4.2) ellipse (0.3cm and 0.3cm);
\end{tikzpicture}}

& \qquad &

\scalemath{0.7}{
\begin{tikzpicture}
  [scale=1,thick,auto=left, every node/.style={circle, inner sep=2pt}] 
  \node[fill=black,label=below:{\Large$E_0V$}] (b0) at (0,0) {};
  \node[fill=black,label=below:{}] (b1) at (1,0)  {};
  \node[fill=black,label=below:{}] (b2) at (2,0)  {};
  \node[fill=black,label=below:{}] (b3) at (3,0) {};
  \node[fill=black,label=below:{}] (b4) at (4,0)  {};
  \node[fill=black,label=below:{}] (b5) at (5,0)  {};
  \node[fill=black,label=below:{}] (b6) at (6,0)  {};
  \node[fill=black,label=below:{\Large$E_dV$}] (b7) at (7,0)  {};
  
  \node[fill=black,label=above:{\Large$E_0^*V$}] (t0) at (0,7)  {};
  \node[fill=black,label=above:{}] (t1) at (1,7)  {};
  \node[fill=black,label=above:{}] (t2) at (2,7)  {};
  \node[fill=black,label=above:{}] (t3) at (3,7)  {};
  \node[fill=black,label=above:{}] (t4) at (4,7)  {};
  \node[fill=black,label=above:{}] (t5) at (5,7)  {}; 
  \node[fill=black,label=above:{}] (t6) at (6,7)  {}; 
  \node[fill=black,label=above:{\Large$E_d^*V$}] (t7) at (7,7)  {}; 

  \node[fill=black,label=left:{}] (l1) at (0,1)  {};
  \node[fill=black,label=left:{}] (l2) at (0,2)  {};
  \node[fill=black,label=left:{}] (l3) at (0,3)  {};
  \node[fill=black,label=left:{}] (l4) at (0,4)  {};
  \node[fill=black,label=left:{}] (l5) at (0,5)  {};
  \node[fill=black,label=left:{}] (l6) at (0,6)  {};

  \node[fill=black,label=right:{}] (r1) at (7,1)  {};
  \node[fill=black,label=right:{}] (r2) at (7,2)  {};
  \node[fill=black,label=right:{}] (r3) at (7,3)  {};
  \node[fill=black,label=right:{}] (r4) at (7,4)  {};
  \node[fill=black,label=right:{}] (r5) at (7,5)  {};
  \node[fill=black,label=right:{}] (r6) at (7,6)  {};

  \node[fill=black,label=right:{}] (d1) at (1,1)  {};
  \node[fill=black,label=right:{}] (d2) at (1.9,1.9)  {};
  \node[fill=black,label=right:{}] (d3) at (2.8,2.8)  {};
  \node[fill=black,label=right:{}] (d4) at (4.2,4.2)  {};
  \node[fill=black,label=right:{}] (d5) at (5.1,5.1)  {};
  \node[fill=black,label=right:{}] (d6) at (6,6)  {};

  \node[fill=black,label=right:{}] (e1) at (1,6)  {};
  \node[fill=black,label=right:{}] (e2) at (1.9,5.1)  {};
  \node[fill=black,label=right:{}] (e3) at (2.8,4.2)  {};
  \node[fill=black,label=right:{}] (e4) at (4.2,2.8)  {};
  \node[fill=black,label=right:{}] (e5) at (5.1,1.9)  {};
  \node[fill=black,label=right:{}] (e6) at (6,1)  {};

  \foreach \from/\to in
  {b0/b1, b1/b2, b2/b3, b3/b4, b4/b5, b5/b6, b6/b7,
  t0/t1, t1/t2, t2/t3, t3/t4, t4/t5, t5/t6, t6/t7,                  
  b0/l1, l1/l2, l2/l3, l3/l4, l4/l5, l5/l6, l6/t0,
  t7/r1, r1/r2, r2/r3, r3/r4, r4/r5, r5/r6, r6/b7,
  t0/e1, e1/e2, e2/e3, e3/e4, e4/e5, e5/e6, e6/b7,
  b0/d1, d1/d2, d2/d3, d3/d4, d4/d5, d5/d6, d6/t7}       
      \draw (\from) -- (\to);

     \draw (4,7) ellipse (0.2cm and 0.2cm);
     \draw (4,7) ellipse (0.3cm and 0.3cm);        
     \draw (5,7) ellipse (0.3cm and 0.3cm);
     
     \draw (4,0) ellipse (0.2cm and 0.2cm);
     \draw (4,0) ellipse (0.3cm and 0.3cm);
     \draw (5,0) ellipse (0.3cm and 0.3cm);
     
     \draw (0,4) ellipse (0.2cm and 0.2cm);
     \draw (0,5) ellipse (0.3cm and 0.3cm);        
     \draw (0,4) ellipse (0.3cm and 0.3cm);
     \draw (0,3) ellipse (0.3cm and 0.3cm);        
     
     \draw (7,4) ellipse (0.2cm and 0.2cm);
     \draw (7,4) ellipse (0.3cm and 0.3cm);        
     
     \draw (1.9,5.1) ellipse (0.2cm and 0.2cm);
     \draw (1.9,5.1) ellipse (0.3cm and 0.3cm);        
     \draw (2.8,4.2) ellipse (0.3cm and 0.3cm);
     
     \draw (1.9,1.9) ellipse (0.2cm and 0.2cm);
     \draw (1.9,1.9) ellipse (0.3cm and 0.3cm);        
     \draw (2.8,2.8) ellipse (0.3cm and 0.3cm);
\end{tikzpicture}}

 \\ 

Action of $\psi_{k+1}$

& \qquad \qquad &

Action of $\tilde\psi_{k+1}$

\end{tabular}
\end{center}

\noindent
In the left diagram, a small circle marks a component $U_i$ on which $\psi_{k+1}$ acts.
The large circles mark the components whose sum contains $\psi_{k+1}U_i$.
The right diagram similarly describes the action of $\tilde\psi_{k+1}$.
\end{minipage}

\smallskip
Recall the basis $W_{-k}, W_{k+1}, \tilde{G}_{k+1}$ $(k \in \mathbb{N})$ from \eqref{eq:3rd basis O}.
We now describe the actions of these basis elements on the decomposition $\{U_i\}_{i=0}^d$.

\begin{corollary}
For $k \in \mathbb{N}$ and $0 \leq i \leq d$, the actions of $W_{-k}$, $W_{k+1}$, $\tilde{G}_{k+1}$ on $U_i$ are described in the table below.
\begin{equation*}
\hspace{-0.3cm}
\renewcommand{\arraystretch}{1.2}
\setlength{\tabcolsep}{4pt}
\begin{tabular}{c|c|c|c}

	decomposition
	& action of $W_{-k}$ on $U_i$
	& action of $W_{k+1}$ on $U_i$ 
	& action of $\tilde{G}_{k+1}$ on $U_i$ \\

	\hline \hline

	$[0,d]$
	& $W_{-k}U_i\subseteq U_i$
	& $W_{k+1}U_i\subseteq U_{i-1}+U_i+U_{i+1}$ 
	& $\tilde{G}_{k+1}U_i\subseteq U_{i-1}+U_i+U_{i+1}$ \\

	$[0^*,d^*]$
	& $W_{-k}U_i\subseteq U_{i-1}+U_i+U_{i+1}$
	& $W_{k+1}U_i\subseteq U_i$ 
	& $\tilde{G}_{k+1}U_i\subseteq U_{i-1}+U_i+U_{i+1}$ \\	

	$[0^*,0]$
	& $W_{-k}U_i\subseteq U_i+U_{i+1}$
	& $W_{k+1}U_i\subseteq U_{i-1}+U_i$
	& $\tilde{G}_{k+1}U_i\subseteq U_{i-1}+U_i+U_{i+1}$ \\

	$[0^*,d]$
	& $W_{-k}U_i\subseteq U_i+U_{i+1}$
	& $W_{k+1}U_i\subseteq U_{i-1}+U_i$
	& $\tilde{G}_{k+1}U_i\subseteq U_{i-1}+U_i+U_{i+1}$ \\

	$[d^*,0]$
	& $W_{-k}U_i\subseteq U_i+U_{i+1}$
	& $W_{k+1}U_i\subseteq U_{i-1}+U_i$
	& $\tilde{G}_{k+1}U_i\subseteq U_{i-1}+U_i+U_{i+1}$ \\

	$[d^*,d]$
	& $W_{-k}U_i\subseteq U_i+U_{i+1}$
	& $W_{k+1}U_i\subseteq U_{i-1}+U_i$
	& $\tilde{G}_{k+1}U_i\subseteq U_{i-1}+U_i+U_{i+1}$ \\	

\end{tabular}
\end{equation*}
\end{corollary}
\begin{proof}
By Proposition~\ref{lem:tm 1st 3rd}, Theorem~\ref{thm:action A_k B_k on Ui}, and Theorem~\ref{thm:action psi_k tilde.psi_k on Ui}.
\end{proof}

\section{Further research}\label{Sec:future-research}

In this section, we give three problems concerning the elements of $\mathcal{O}$ displayed in \eqref{eq:module elements}.
For convenience, we recall these elements:
\begin{equation}\label{eq:problem elements}
	A_k, \qquad B_k, \qquad \psi_{k+1}, \qquad \tilde\psi_{k+1}, \qquad k\in\mathbb{N}.
\end{equation}
Throughout this section, let $V$ denote a finite-dimensional irreducible $\mathcal{O}$-module.

\begin{problem}\label{prob1}
Determine which of the elements in \eqref{eq:problem elements} are diagonalizable on $V$.
\end{problem}

\begin{problem}\label{prob2}
This problem has four parts.
\begin{itemize}
	\item[(i)]
	Assume that $A_k$ is diagonalizable on $V$ for every $k\in\mathbb{N}$.
	Find an explicit basis for $V$ consisting of common eigenvectors for the $\{A_k\}_{k \in \mathbb{N}}$.

	\item[(ii)]
	Assume that $B_k$ is diagonalizable on $V$ for every $k\in\mathbb{N}$.
	Find an explicit basis for $V$ consisting of common eigenvectors for the $\{B_k\}_{k \in \mathbb{N}}$.

	\item[(iii)]
	Assume that $\psi_{k+1}$ is diagonalizable on $V$ for every $k\in\mathbb{N}$.
	Find an explicit basis for $V$ consisting of common eigenvectors for the $\{\psi_{k+1}\}_{k \in \mathbb{N}}$. 

	\item[(iv)]
	Assume that $\tilde\psi_{k+1}$ is diagonalizable on $V$ for every $k\in\mathbb{N}$.
	Find an explicit basis for $V$ consisting of common eigenvectors for the $\{\tilde\psi_{k+1}\}_{k \in \mathbb{N}}$.
\end{itemize}
\end{problem}

\begin{problem}\label{prob3}
Assume that every element in \eqref{eq:problem elements} is diagonalizable on $V$.
Each {\rm(i)}--{\rm(iv)} in Problem~\ref{prob2} yields an explicit eigenbasis for $V$.
Determine the transition matrices between these four bases.
\end{problem}



\section*{Acknowledgement}
The author would like to express his deep gratitude to Paul Terwilliger for many valuable discussions and insightful suggestions during the preparation of this paper.


\end{document}